\documentclass{amsart}

\usepackage{amssymb}

\allowdisplaybreaks

\newtheorem{theorem}{Theorem}[section]
\newtheorem{lemma}[theorem]{Lemma}
\newtheorem{proposition}[theorem]{Proposition}
\newtheorem{corollary}[theorem]{Corollary}

\theoremstyle{definition}

\theoremstyle{remark}
\newtheorem{remark}[theorem]{Remark}

\newcommand{\Z}{\mathbb{Z}}
\newcommand{\R}{\mathbb{R}}
\newcommand{\C}{\mathbb{C}}

\newcommand{\Pol}{\mathbb{P}}

\newcommand{\Po}{\mathrm{P}}
\newcommand{\Su}{\mathrm{S}}
\newcommand{\Di}{\mathrm{D}}
\newcommand{\OH}{\mathrm{L}}

\begin{document}

\title[Classical hypergeometric
polynomials] {A new approach to the theory of classical
hypergeometric polynomials}

\author[Marco and Parcet]
{Jos\'{e} Manuel Marco and Javier Parcet}

\address{Department of Mathematics, Universidad Aut\'{o}noma de
Madrid, Madrid 28049, Spain}

\email{javier.parcet@uam.es}

\subjclass{Primary 33D15, 33D45}

\date{}

\keywords{P-sequence, hypergeometric operator, Taylor and
Rodrigues formula}

\maketitle

\begin{abstract}
In this paper we present a unified approach to the spectral
analysis of an hypergeometric type operator whose eigenfunctions
include the classical orthogonal polynomials. We write the
eigenfunctions of this operator by means of a new Taylor formula
for operators of Askey-Wilson type. This gives rise to some
expressions for the eigenfunctions, which are unknown in such a
general setting. Our methods also give a general Rodrigues formula
from which several well known formulas of Rodrigues type can be
obtained directly. Moreover, other new Rodrigues type formulas
come out when seeking for regular solutions of the associated
functional equations. The main difference here is that, in
contrast with the formulas appearing in the literature, we get
non-ramified solutions which are useful for applications in
combinatorics. Another fact, that becomes clear in this paper, is
the role played by the theory of elliptic functions in the
connection between ramified and non-ramified solutions.
\end{abstract}

\section*{Introduction}

A large number of discrete problems in graph theory, group
representation theory, discrete harmonic analysis or finite
symmetric spaces involve the spectral analysis of certain
difference operators. Namely, given a positive integer $d$ and
denoting by $\mathrm{V}_d$ the space of functions $f: \{0,1,
\ldots, d\} \rightarrow \C$, we are interested in the
eigenfunctions of linear operators $\Lambda: \mathrm{V}_d
\rightarrow \mathrm{V}_d$ of the following type $$\Lambda f(t) =
\sigma_+(t) (f(t+1) - f(t)) - \sigma_-(t) (f(t) - f(t-1))$$ where
$\sigma_{\pm} \in \mathrm{V}_d$ and $\sigma_-(0) = \sigma_+(d) = 0$.
For instance, the work \cite{MP1} illustrates how this kind of
eigenvalue problems become useful to solve some a priori unrelated
combinatorial problems. Let us write $\Pol_k[x]$ for the space of
complex polynomials of degree $\le k$ in one variable. Explicit
solutions of this spectral analysis are available when there
exists an injective mapping $x: t \in \{0,1, \ldots, d\} \mapsto
x_t \in \C$ such that the following conditions hold
\begin{itemize}
\item[(a)] We have $x_{t+1} + x_{t-1} = f(x_t)$ for some $f \in
\Pol_1[x]$ and $1 \le t \le d-1$.
\item[(b)] The subspaces $\mathrm{V}_{d,k} = \{f \circ x:
\, f \in \Pol_k[x]\}$ are $\Lambda$-invariant for $0 \le k \le
d-1$.
\end{itemize}
By a well-known result of Leonard \cite{L} --see also the
reference \cite{BCN}-- we have $(\textrm{b} \Rightarrow
\textrm{a})$ whenever $d > 7$. When this is the case, the
eigenfunctions of $\Lambda$ can be obtained by applying the theory
of basic hypergeometric series. Alternatively, Nikiforov, Suslov
and Uvarov study in \cite{NSU} the spectral decomposition of
$\Lambda$ using properties (a) and (b) and some specific methods.
This gives rise to a whole new approach to the theory of classical
orthogonal polynomials. Following any of these methods, two very
relevant expressions for the eigenfunctions of $\Lambda$ can be
obtained. First, as we have recalled above, we can write the
eigenfunctions in terms of $q$-hypergeometric polynomials. Second,
the eigenfunctions can be written in general by means of a
Rodrigues type formula. If we look at the approach which uses the
theory of basic hypergeometric series, a long list of subcases
must be checked. For instance, the reader is referred to the book
of Bannai and Ito \cite{BI}. Besides, in that reference the
authors do not provide any Rodrigues type formula for a problem
which in essence is equivalent to the one considered here. On the
other hand, in the approach of \cite{NSU}, the derivation of
Rodrigues type formulas is cumbersome. Moreover, at the very final
stages, the authors need to invoke the Watson-Sears transformation
and other formulas of the theory of $q$-series to get the
hypergeometric expressions.

In this paper we propose an intermediate approach which, from our
point of view, is quite simpler and provides a somehow natural
classification of the whole problem which seems to be unknown.
Throughout the paper, we shall also obtain some results and
remarks which, as far as we know, are new in the theory. Following
Nikiforov, Suslov and Uvarov, we shall study the eigenfunctions of
$\Lambda$ with the aid of properties (a) and (b), so that the
theory of basic hypergeometric series does not play a central
role. However, in contrast with \cite{NSU}, we reformulate the
operator $\Lambda$ in terms of operators of Askey-Wilson type, so
that we work with a generalized hypergeometric operator $\OH$. By
an Askey-Wilson operator $\Di$, we do not only mean the operators
introduced by Richard Askey and James A. Wilson in \cite{AW}.
Actually, we shall also be working with other kind of operators
strongly related with those, such as the one introduced by Hahn.
Besides, we shall use the companion operator $\Su$ in the sense of
Magnus \cite{M}. In particular, the hypergeometric operator we
shall work with has the form $$\OH = \sigma \Di^2 + \tau \Su
\Di,$$ where $\sigma \in \Pol_2[x]$ and $\tau \in \Pol_1[x]$. The
operators $\Su$ and $\Di$ are defined by means of a two variable
quadratic polynomial $\Po$ which in our case will be symmetric. In
our context, a translation of the difference operator $\Delta f(t)
= f(t+1) - f(t)$ and the classical derivative will be particular
cases. These Askey-Wilson type operators act over polynomials
--also over analytic functions-- without any explicit reference to
a particular sequence $x_t$ satisfying properties (a) and (b).
Formulated in these terms, it turns out that the operator
$\Lambda$ becomes more intrinsic and, as we shall see, this
provides some relevant advantages.

In order to study the eigenvalue problem posed above, we
previously state some results for the Askey-Wilson operators. We
first give seven canonical forms under affine transformations
which allow us to provide a very natural classification of the
Askey-Wilson operators. We also develop the Leibniz and
commutation formulas for the Askey-Wilson operator and its
companion, the corresponding Taylor formula and Cauchy's integral
formula. Results of this kind are interesting on their own right
and can be used to obtain new proofs of well-known summation
formulas for some basic hypergeometric series. We have not
included here these applications since they will be the subject of
a forthcoming publication \cite{MP2}. The use of these techniques
and mainly our Taylor formula, allow us to obtain simple
expressions for the eigenfunctions of the hypergeometric operator
$\OH$ in a very general context, see Theorem
\ref{Teorema-Expresiones-Polinomicas}. In particular, we do not
need to invoke any Rodrigues type formula or any result from the
theory of basic hypergeometric polynomials, such as the
Watson-Sears transformations used in \cite{NSU}. Moreover, the
identities given in Theorem \ref{Teorema-Expresiones-Polinomicas}
seem to be unknown in such a general setting and, as we shall see,
can be easily applied to particular cases. The expressions of the
eigenfunctions in terms of $q$-hypergeometric functions arise from
this result. On the other hand, Theorem \ref{Teorema-Rodrigues}
provides an algorithm to obtain Rodrigues type formulas which
depend on the resolution of certain functional equation. Our
methods in Theorem \ref{Teorema-Rodrigues} are closely related to
those used in \cite{NSU}. However, the final result of Nikiforov,
Suslov and Uvarov is not so clear since they do not provide the
relation with the Askey-Wilson operators. Moreover, it is not easy
to see the dependence of their Rodrigues formulas with respect to
the function $x_t$ which appears in properties (a) and (b). Our
Rodrigues formulas are also related to the ones obtained by Askey
and Wilson in \cite{AW} for the Askey-Wilson polynomials.

The canonical Askey-Wilson operators mentioned above lead us to
consider five canonical families of eigenvalue problems. For some
reasons that will become clear throughout the paper, we call these
families as follows: continuous, arithmetic, quadratic, geometric
and trigonometric. The corresponding eigenfunctions include as
generic cases the Jacobi, Hahn, Wilson, $q$-Hahn and Askey-Wilson
polynomials respectively. We shall parameterize each family by
means of one or two auxiliary polynomials. In the trigonometric
case, we get a Laurent polynomial. Using this, we rewrite the
expressions of Theorem \ref{Teorema-Expresiones-Polinomicas} in
terms of the basic hypergeometric functions. Besides, we will
study the functional equation of Theorem \ref{Teorema-Rodrigues}
to obtain explicit Rodrigues formulas for each canonical case. Our
approach looks mainly at the regularity of the solutions of the
functional equation. This gives rise to new Rodrigues type
formulas which do not appear in the literature. Quite
surprisingly, apart from the continuous case, we can always obtain
non-ramified solutions. More specifically, we get either
meromorphic solutions or solutions with one essential singularity.
In general, there is not uniqueness of solutions, so that we
investigate a great variety of solutions. Although we shall only
give explicit identities for some relevant generic cases, it will
become clear that our methods can be applied in all cases.
Moreover, in contrast with \cite{NSU}, our methods allow us to
obtain such relations directly. That is, we do not need to take
limits. We also point out that our Rodrigues formula for the
Askey-Wilson polynomials is different from the one given by Askey
and Wilson in \cite{AW}. From our point of view, looking at the
regularity of the functions appearing in Rodrigues formula, Askey
and Wilson obtained ramified solutions while we get meromorphic
ones. By using the non-uniqueness of the associated functional
equation, we show how our Rodrigues formula can be used to get the
orthogonality relations given in \cite{AW}. In fact, the
techniques we develop here can be used to obtain discrete and
absolutely continuous measures with respect to which, the
corresponding eigenfunctions become pairwise orthogonal. Besides,
one of the most relevant features of our techniques is that the
role of the theory of elliptic functions in the problem becomes
quite clear. More concretely, our approach divides the analysis of
Rodrigues type formulas into two different parts. First, we seek
for regular Rodrigues type formulas avoiding ramified functions.
Formulas of this kind are useful for applications in
combinatorics. Second, with the aid of the theory of elliptic
functions, we explore the non-uniqueness of solutions to obtain
the original Rodrigues type formulas given in \cite{AW}.

Throughout this paper, we use the customary notation \cite{GR} for
$q$-shifted factorials and basic hypergeometric series. We also
use an auxiliary non-zero parameter $\lambda$ with $q =
\lambda^2$. The restriction $|q| < 1$ is only assumed in some
explicit expressions. For instance, the main Theorems
\ref{Teorema-Expresiones-Polinomicas} and \ref{Teorema-Rodrigues}
hold without that restriction.

\section{$\Po$-sequences}
\label{Section1}

In this section we introduce the notions of $\Po$-sequence and
$\Po$-function associated to a given symmetric polynomial $\Po$.
As we shall see, these notions become very useful in order to give
a somehow natural classification of the family of operators of
Askey-Wilson type. Given $a,b,c$ complex numbers, let us consider
the symmetric polynomial in two variables $\Po(x,y) = x^2 + y^2 -
2axy - 2b(x+y) + c$. For any such polynomial we have $$\Po(x,y) =
y^2 - 2A(x)y + B(x),$$ where $A(x) = ax + b$ and $B(x) = x^2 - 2bx
+ c$. Let $\frac{1}{2} \Z$ stand for the set of integer and
semi-integer numbers. A sequence of complex numbers $(x_t)$, with
the index $t$ running over $\frac{1}{2} \Z$, will be called a
$\Po$-\textbf{sequence} if $$\Po(x_t,y) = (y - x_{t +
\frac{1}{2}}) (y - x_{t - \frac{1}{2}}) \qquad \mbox{for all $t
\in \frac{1}{2} \Z$}.$$ If $(x_t)$ is a $\Po$-sequence and $t_0
\in \frac{1}{2} \Z$, the translated sequence $(y_t)$ with $y_t =
x_{t - t_0}$ is also a $\Po$-sequence. We shall say that $x_0$ is
the \textbf{base point} of the $\Po$-sequence $(x_t)$. Given a
complex number $\xi$, a $\Po$-sequence with base point $\xi$ can
be constructed recursively. Moreover, this $\Po$-sequence is
unique once we have decided which root of $\Po(\xi,y)$ is
$x_{1/2}$ and which one is $x_{-1/2}$. In particular, we obtain
two $\Po$-sequences with base point $\xi$ which coincide when
$\Po(\xi,y)$ has a double root. If $(x_t)$ is such a
$\Po$-sequence, the other $\Po$-sequence with base point $\xi$ is
$(y_t)$ where $y_t = x_{-t}$.

\begin{remark}
Obviously, the symmetry of $\Po$ plays an essential role in the
iteration.
\end{remark}

As we have seen, $\Po$ can be regarded as a polynomial in the
variable $y$ with coefficients depending on $x$. This allows us to
consider the discriminant of $\Po$ as a function of the variable
$x$ which, up to a constant factor, has the form $$ \delta(x) =
A(x)^2 - B(x) = (a^2-1)x^2 + 2b(a+1)x + b^2 - c.$$ By the
definition of $\Po$-sequence, the following identities hold for
any $t \in \frac{1}{2} \Z$
\begin{eqnarray} \label{ABdelta}
\nonumber A(x_t) & = & \textstyle \frac{1}{2} (x_{t + \frac{1}{2}}
+ x_{t - \frac{1}{2}}) \\ B(x_t) & = & x_{t + \frac{1}{2}} x_{t -
\frac{1}{2}} \\ \nonumber \delta(x_t) & = & \textstyle
\frac{1}{4}(x_{t + \frac{1}{2}} - x_{t - \frac{1}{2}})^2.
\end{eqnarray}

\begin{remark}
(\ref{ABdelta}) characterizes $\Po$-sequences with $(A(x_t),
B(x_t))$ or $(A(x_t),\delta(x_t))$.
\end{remark}

Each symmetric polynomial $\Po(x,y) = x^2 + y^2 - 2axy -2b(x+y) +
c$ is completely determined by its coefficients $a,b,c$. Given a
parameterized curve $$\begin{array}{rrcl} \gamma: & \R &
\longrightarrow & \C^3 \\ & s & \longmapsto & (a(s),b(s),c(s)),
\end{array}$$ we consider the polynomials $\Po_s(x,y) = x^2 + y^2
- 2 a(s) xy - 2 b(s) (x+y) + c(s)$. Then, once we have fixed a
complex number $x_0$, we can construct for any $s \in \R$ a
$\Po_s$-sequence $(x_t(s))$ with base point $x_0$ such that the
whole family of $\Po_s$-sequences depend continuously on the
parameter $s$. Namely, if $\delta_s$ stands for the discriminant
of $\Po_s$, we can choose the value of $\sqrt{\delta_s(x_0)}$ in
such a way that we obtain a continuous function on $s$. This
allows us to define $$x_{\pm \frac{1}{2}}(s) = A_s(x_0) \pm
\sqrt{\delta_s(x_0)}$$ and then we construct the $\Po_s$-sequence
recursively. The continuous dependence along curves expressed
above shows that any $\Po$-sequence with $a = \pm 1$ can be
approximated by $\Po$-sequences with $a \neq \pm 1$. This
continuity argument will be very useful in the sequel and we shall
use it in what follows with no further reference.

Let $\mathcal{P}$ be the set of polynomials $\Po(x,y) = x^2 + y^2
- 2 a xy - 2 b (x+y) + c$. The set $\mathcal{P}$ is an affine
subspace of the vector space $\mathbb{SP}_2[x,y]$ of symmetric
polynomials in two complex variables having degree $\le 2$. The
affine group of the complex plane $\mbox{Aff}(\C)$ naturally acts
on $\mathcal{P}$ by the action $\mathcal{A}_1: \mbox{Aff}(\C)
\times \mathcal{P} \rightarrow \mathcal{P}$ given by $$(g \cdot
\Po) (x,y) = \zeta^{2} \Po(g^{-1}(x),g^{-1}(y)),$$ where $g(x) =
\zeta x + \eta$. The parameter $a$, which will be called the
\textbf{main coefficient} of $\Po$, is invariant under the action
of $\mbox{Aff}(\C)$ and the discriminant $\delta$ is transformed
by the following rule
\begin{equation} \label{Ecuacion-Orbita2}
\delta^{g \cdot \Po} (x) = \zeta^2 \delta^{\Po} (g^{-1}(x)).
\end{equation}
The relation above defines another action $\mathcal{A}_2:
\mbox{Aff}(\C) \times \mathbb{P}_2[x] \rightarrow \mathbb{P}_2[x]$
on the vector space $\mathbb{P}_2[x]$ of complex polynomials in
one variable with degree $\le 2$. As we shall see throughout the
paper, all the notions associated to the elements $\Po$ of the set
$\mathcal{P}$ are naturally preserved under the action of
$\mbox{Aff}(\C)$. For instance, if we are given a $\Po$-sequence
$(x_t)$, we have that $(g(x_t))$ is a $(g \cdot \Po)$-sequence. A
$\Po$-sequence $(x_t)$ is an even function of the variable $t$ if
and only if its base point $x_0$ is a root of the discriminant
$\delta$. For this reason we define the \textbf{even points} of
$\Po$ to be the roots of its discriminant $\delta$. The number of
distinct even points of $\Po$ will be denoted by $ev_{\Po}$. Now
we show that the main coefficient of $\Po$ and the parameter
$ev_{\Po}$ are relevant invariants in the orbits of the action
$\mathcal{A}_1$.

\begin{proposition} \label{Proposicion-Orbitas}
Given a pair of symmetric polynomials $\Po_1, \Po_2 \in
\mathcal{P}$ with main coefficients $a_j$ and discriminants
$\delta_j$ for $j = 1,2$, the following are equivalent:
\begin{itemize}
\item[\textnormal{(a)}] The symmetric polynomials $\Po_1$ and
$\Po_2$ belong to the same $\mathcal{A}_1$-orbit.
\item[\textnormal{(b)}] The discriminants $\delta_1$ and
$\delta_2$ have the same degree, $ev_{\Po_1} = ev_{\Po_2}$ and
$a_1 = a_2$.
\end{itemize}
\end{proposition}

\begin{proof}
The implication $(\mbox{a} \Rightarrow \mbox{b})$ follows from
(\ref{Ecuacion-Orbita2}) and the recalled fact that the main
coefficient remains invariant under the action of
$\mbox{Aff}(\C)$. Let us see that $(\mbox{b} \Rightarrow
\mbox{a})$. If (b) holds, it is clear that $\delta_1$ and
$\delta_2$ belong to the same $\mathcal{A}_2$-orbit. In
particular, there exists $g_{12} \in \mbox{Aff}(\C)$ such that the
polynomial $\Po = g_{12} \cdot \Po_1$ has discriminant $\delta_2$.
If $a \neq -1$ this implies $\Po = \Po_2$ and we are done.
Finally, for $a = -1$ we can always choose $g_1, g_2 \in
\mbox{Aff}(\C)$ such that $$(g_j \cdot \Po_j)(x,y) = x^2 + y^2 +
2xy + c_j \qquad \mbox{for} \qquad j = 1,2$$ where $c_j = 1$ if
$\delta_j \neq 0$ and $c_j = 0$ otherwise. Since $\delta_1 = 0$ if
and only if $\delta_2 = 0$ by hypothesis, $\Po_1$ and $\Po_2$
belong to the same $\mathcal{A}_1$-orbit and the proof is
completed.
\end{proof}

Using Proposition \ref{Proposicion-Orbitas} we can write any
polynomial $\Po$ in $\mathcal{P}$, and its corresponding
$\Po$-sequences, in a \textbf{canonical form}. Namely, in the
following table we write all the $\mathcal{A}_1$-orbits by
considering all the possible combinations of
$\mathcal{A}_1$-invariants: the main coefficient $a$, the degree
$dg(\delta)$ of the discriminant and the number of distinct even
points $ev_{\Po}$. We have chosen a canonical symmetric polynomial
$\Po$ in each $\mathcal{A}_1$-orbit. Then, we obtain the
corresponding $\Po$-sequences. Our election for the canonical
forms has given priority to the simplicity of the resulting
$\Po$-sequence and it is somehow arbitrary. Finally, recalling
that $\mbox{Aff}(\C)$ can be regarded as the semi-direct product
$\mathfrak{T}(\C) \rtimes GL(\C)$ where $\mathfrak{T}(\C)$ is the
group of translations in the complex plane, we include the
isotropy subgroup $\mathbf{I}_{\Po}$ associated to each orbit. In
the following table, we use a new parameter $\lambda \neq 0$
defined by the relation $\lambda = a \pm \sqrt{a^2 -1}$. In other
words, $$a = \frac{\lambda + \lambda}{2}^{-1}.$$ Obviously, the
product of the two values assigned to $\lambda$ is 1. Moreover, we
have $$\begin{array}{rcl} \lambda = 1 & \Leftrightarrow & a = 1 \\
\lambda = -1 & \Leftrightarrow & a = -1
\end{array}$$ As we have noticed above, these values of $\lambda$ can be
regarded as limiting cases.

\begin{center}
\begin{tabular}{|r|c|c|c|c|c|c|} \multicolumn{1}{c}{} &
\multicolumn{1}{c}{\scriptsize \textbf{\emph{a}}} &
\multicolumn{1}{c}{\scriptsize
\textbf{\emph{$\mbox{dg}(\delta)$}}} &
\multicolumn{1}{c}{\scriptsize \textbf{\emph{$\mbox{ev}_{\Po}$}}}
& \multicolumn{1}{c}{\scriptsize \textbf{Canonical form}} &
\multicolumn{1}{c}{\scriptsize \textbf{$\Po$-sequence}} &
\multicolumn{1}{c}{\scriptsize \textbf{$\textrm{I}_{\Po}$}}
\\ \cline{1-7} \scriptsize \textbf{T} & \scriptsize $\neq \pm 1$ &
\scriptsize $2$ & \scriptsize $2$ & \scriptsize $x^2 + y^2 - 2 a
xy + a^2 - 1$ & \scriptsize $\frac{1}{2} (\lambda^{2t} u +
\lambda^{-2t} u^{-1})$ & \scriptsize $\{\pm 1\}$ \\ \cline{1-7}
\scriptsize \textbf{G} & \scriptsize $\neq \pm 1$ & \scriptsize
$2$ & \scriptsize $1$ & \scriptsize $x^2 + y^2 - 2 a xy$ &
\scriptsize $\lambda^{\pm 2t} x_0$ & \scriptsize $GL(\C)$ \\
\cline{1-7} \scriptsize \textbf{Q} & \scriptsize $1$ & \scriptsize
$1$ & \scriptsize $1$ & \scriptsize $(x-y)^2 - \frac{1}{2}(x+y) +
\frac{1}{16}$ & \scriptsize $t^2 + 2tu + u^2$ & \scriptsize $1$
\\ \cline{1-7} \scriptsize \textbf{A} & \scriptsize $1$ &
\scriptsize $0$ & \scriptsize $0$ & \scriptsize $(x-y)^2 -
\frac{1}{4}$ & \scriptsize $x_0 \pm t$ & \scriptsize
$\mathfrak{T}(\C) \rtimes \{\pm 1\}$ \\ \cline{1-7} \scriptsize
\textbf{C} & \scriptsize $1$ & \scriptsize $- \infty$ &
\scriptsize $\infty$ & \scriptsize $(x-y)^2$ & \scriptsize $x_0$ &
\scriptsize $\mbox{Aff}(\C)$ \\ \cline{1-7} \scriptsize \textbf{O}
& \scriptsize $-1$ & \scriptsize $0$ & \scriptsize $0$ &
\scriptsize $(x+y)^2 - \frac{1}{4}$ & \scriptsize $(-1)^{2t}(x_0
\pm t)$ & \scriptsize $\{\pm 1\}$ \\ \cline{1-7} \scriptsize
\textbf{E} & \scriptsize $-1$ & \scriptsize $- \infty$ &
\scriptsize $\infty$ & \scriptsize $(x+y)^2$ & \scriptsize
$(-1)^{2t}x_0$ & \scriptsize $GL(\C)$ \\ \cline{1-7}
\end{tabular}
\\ \vskip10pt
\textsc{Table I.} Canonical $\Po$-sequences. \vskip6pt
\end{center}

The first two canonical forms include infinitely many
$\mathcal{A}_1$-orbits, one for each value of the main coefficient
$a$. There are only seven inequivalent canonical forms. The
capital letters in the left column are acronyms of the names we
have adopted for the canonical forms. These names have been chosen
attending at the shape of the canonical $\Po$-sequences:
\textbf{Trigonometric}, \textbf{Geometric}, \textbf{Quadratic},
\textbf{Arithmetic}, \textbf{Constant}, \textbf{Oscillating} and
\textbf{Even}. For some reasons that will become clear later, we
shall also refer to the Constant form as the \textbf{Continuous}
canonical form.

\begin{remark}
It is clear that the polynomials belonging to the
$\mathcal{A}_1$-orbits associated to the trigonometric canonical
form are dense in $\mathcal{P}$. Therefore, the trigonometric
canonical form is the most relevant one since any other canonical
form can be regarded as a limit of families of polynomials in this
particular orbit.
\end{remark}

\begin{remark} \label{Observacion-P-sucesiones}
Taking affine transformations from the canonical $\Po$-sequences,
it is not difficult to check that every $\Po$-sequence $(x_t)$ has
one of the following forms:
\begin{itemize}
\item If $a \neq \pm 1$, then we have $x_t =
k_0 + k_1 \lambda^{2t} + k_2 \lambda^{-2t}$ where
\begin{eqnarray*}
x_0 & = & k_0 + k_1 + k_2 \\ b & = & (1 - a) k_0 \\ c & = & 2 b
k_0 + 4 k_1 k_2 (a^2 - 1).
\end{eqnarray*}

\item If $a = 1$, then we have $x_t = k_0 + k_1 t + k_2 t^2$ where
\begin{eqnarray*}
x_0 & = & k_0 \\ b & = & \textstyle \frac{1}{4} k_2 \\ c & = &
\textstyle \frac{1}{16} k_2^2 - \frac{1}{4} k_1^2 + k_0 k_2. \quad
{}
\end{eqnarray*}

\item If $a = -1$, then we have $x_t = k_0 + (-1)^{2t} (k_1 +
k_2 t)$ where
\begin{eqnarray*}
x_0 & = & k_0 + k_1 \\ b & = & 2 k_0 \\ c & = & \textstyle 4 k_0^2
- \frac{1}{4} k_2^2. \qquad \qquad \ \ {}
\end{eqnarray*}
\end{itemize}
These relations already appeared in the book \cite{NSU} of
Nikiforov, Suslov and Uvarov.
\end{remark}

\begin{proposition} \label{Proposicion-Formulario}
Every $\Po$-sequence $(x_t)$ satisfies the following relations,
with the obvious limits for $\lambda = \pm 1$, for any $s \in
\frac{1}{2} \Z$
\begin{eqnarray*}
\frac{x_{t+s} + x_{t-s}}{2} & = & \frac{\lambda^{2s} +
\lambda^{-2s}}{2} \ x_t + \frac{\lambda^{2s} + \lambda^{-2s} -
2}{\lambda + \lambda^{-1} - 2} \ b \\ x_{t+s} - x_{t-s} & = &
\frac{\lambda^{2s} - \lambda^{-2s}}{\lambda - \lambda^{-1}} \
(x_{t + \frac{1}{2}} - x_{t - \frac{1}{2}}).
\end{eqnarray*}
\end{proposition}

\begin{proof}
By continuity, it suffices to consider $\Po$-sequences with $a
\neq \pm 1$. Moreover, by a translation in the variable $t$, we
can assume that $t = 0$. Finally, under these new hypothesis, the
result follows easily from Remark \ref{Observacion-P-sucesiones}.
\end{proof}

\begin{corollary} \label{Corolario-Formulario}
Every $\Po$-sequence $(x_t)$ satisfies the following relations:
\begin{itemize}
\item[\textnormal{(a)}] $\displaystyle x_{t+1} - x_{t-1} = 2 a \,
(x_{t + \frac{1}{2}} - x_{t - \frac{1}{2}})$.
\item[\textnormal{(b)}] $\displaystyle 8 a \, \delta(x_t) = (x_{t+1}
- x_{t-1}) (x_{t + \frac{1}{2}} - x_{t - \frac{1}{2}})$.
\item[\textnormal{(c)}] $\displaystyle x_{t+1} - 2 x_t + x_{t-1} =
4 (a+1) (A(x_t) - x_t)$.
\end{itemize}
\end{corollary}

Given a symmetric polynomial $\Po \in \mathcal{P}$, an entire
function $\theta: \C \rightarrow \C$ will be called a
\textbf{$\Po$-function} if the sequence $$\ldots \theta(z - 1),
\theta(z - \textstyle \frac{1}{2}), \theta(z), \theta(z +
\frac{1}{2}), \theta(z + 1), \ldots$$ is a $\Po$-sequence for any
$z \in \C$. To be consistent with our notation, we shall write
$\theta_t$ for $\theta(t)$. $\Po$-functions and $\Po$-sequences
can be treated in the same fashion. For instance, an entire
function $\theta: \C \rightarrow \C$ is a $\Po$-function if and
only if the following relations hold for any $t \in \C$
$$A(\theta_t) = \textstyle \frac{1}{2} (\theta_{t + \frac{1}{2}} +
\theta_{t - \frac{1}{2}}) \qquad \mbox{and} \qquad B(\theta_t) =
\theta_{t + \frac{1}{2}} \theta_{t - \frac{1}{2}}.$$

\begin{remark}
In what follows, when dealing with $\Po$-functions, we shall write
$\lambda^t$ as an abbreviation of $e^{t \omega}$ for some election
$\omega$ of the logarithm $\log \lambda$ and any complex number
$t$. We shall also make this slight abuse of notation for $q^t$
where $q = \lambda^2$.
\end{remark}

\section{Taylor formula for the Askey-Wilson operator}
\label{Section2}

In this section we develop a Taylor formula for operators of
Askey-Wilson type, which will we shall apply in Sections
\ref{Section3} and \ref{Section4} to give explicit expressions for
the polynomic eigenfunctions of the associated hypergeometric
operator.

\subsection{Operators $\Su$ and $\Di$}

If $\Pol[x]$ stands for the space of complex polynomials in one
variable, let us denote by $\Pol_k[x]$ the subspace of polynomials
of degree $\le k$. Given $f \in \Pol[x]$, we define $f_d(x,y)$ as
the polynomial which coincides with $$\frac{f(x) - f(y)}{x-y}$$
whenever $x \neq y$. We also introduce the polynomial $f_s(x,y) =
\frac{1}{2} (f(x) + f(y))$. Both $f_s(x,y)$ and $f_d(x,y)$ are
symmetric polynomials in two variables. On the other hand, given a
symmetric polynomial $\Po(x,y) = x^2 + y^2 - 2 a xy - 2 b (x+y) +
c$, we shall denote by $u(x)$ and $v(x)$ the roots of $\Po(x,y)$
in the variable $y$. That is, we have $\Po(x,y) = (y - u(x)) (y -
v(x))$ for all $x \in \C$. Then we define the \textbf{operators}
$\Su$ and $\Di$ as follows $$\Su f(x) = f_s(u(x),v(x)) \qquad
\mbox{and} \qquad \Di f(x) = f_d(u(x),v(x)).$$ As it is recalled
by Magnus \cite{M}, the divided difference operator $\Di$ is
essentially the Askey-Wilson operator while $\Su$ will be called
the companion operator to $\Di$. Any symmetric polynomial in two
variables provides a polynomial in the variable $x$ when it is
evaluated at $(u(x),v(x))$. In particular, we have that $\Su f \in
\Pol[x]$ and $\Di f \in \Pol[x]$ for any $f \in \Pol[x]$. That is,
$\Su, \Di \in \mbox{End}(\Pol[x])$ are linear mappings in
$\Pol[x]$.

\begin{remark} \label{Observacion-SD-Accion}
$\Su$ and $\Di$ are naturally transformed under the action of
$\mbox{Aff}(\C)$. That is, if $\Su^{\Po}$ and $\Di^{\Po}$ stand
for the operators associated to the symmetric polynomial $\Po$ and
we take some $g = (\eta, \zeta) \in \mathfrak{T}(\C) \rtimes
GL(\C) = \mbox{Aff}(\C)$, then we have
$$g \ \Su^{\Po} g^{-1} = \Su^{g \cdot \Po} \qquad \mbox{and}
\qquad g \ \Di^{\Po} g^{-1} = \zeta \, \Di^{g \cdot \Po}$$ where
we identify $g \in \mbox{Aff}(\C)$ with the induced mapping $(g
\cdot f) (x) = f (g^{-1} (x))$.
\end{remark}

\begin{proposition} \label{Proposicion-TerminoDominante}
We have $\Su(\Pol_k[x]) \subset \Pol_k[x]$ and $\Di(\Pol_k[x])
\subset \Pol_{k-1}[x]$ for $k \ge 1$. Moreover, if $f(x) = x^k$ we
obtain the following expressions
\begin{eqnarray*}
\Su f(x) & = & \frac{\lambda^k + \lambda^{-k}}{2} \ x^k + \ldots
\\ \Di f(x) & = & \frac{\lambda^k - \lambda^{-k}}{\lambda -
\lambda^{-1}} \ x^{k-1} + \ldots
\end{eqnarray*}
where the dots stand for terms of lower degree. If $f \in
\Pol_0[x]$, $\Su f = f$ and $\Di f = 0$.
\end{proposition}

\begin{proof}
The action of $\Su$ and $\Di$ on $\Pol_0[x]$ is trivial. Now, let
$f(x) = x^k$ for some $k \ge 1$. Then, since the roots of
$\Po(x,y)$ are given by $A(x) \pm \sqrt{\delta(x)}$, we have $$\Di
f(x) = \frac{\big( A(x) + \sqrt{\delta(x)} \big)^k - \big( A(x) -
\sqrt{\delta(x)} \big)^k}{2 \sqrt{\delta(x)}}$$ whenever
$\delta(x) \neq 0$. The binomial theorem gives
\begin{eqnarray*}
\Di f(x) & = & \frac{1}{2} \sum_{j=0}^k \binom{k}{j} A(x)^{k-j}
(\sqrt{\delta(x)})^{j-1} (1 - (-1)^j) \\ & = & \frac{x^{k-1}}{2
\sqrt{a^2 - 1}} \sum_{j=0}^k \binom{k}{j} a^{k-j} (\sqrt{a^2 -
1})^j (1 - (-1)^j) + \ldots \\ & = & \frac{x^{k-1}}{2 \sqrt{a^2 -
1}} \big[ (a + \sqrt{a^2 - 1})^k - (a - \sqrt{a^2 - 1})^k \big] +
\ldots
\end{eqnarray*}
where the dots stand for terms of lower degree. This new
expression is now valid for any $x \in \C$. In particular, we have
obtained $$\Di f(x) = \frac{\lambda^k - \lambda^{-k}}{\lambda -
\lambda^{-1}} \ x^{k-1} + \ldots$$ Similar arguments provide the
expression given above for $\Su f(x)$.
\end{proof}

\begin{remark}
Both operators $\Su$ and $\Di$ depend polynomially on $\Po$. In
particular, Proposition \ref{Proposicion-TerminoDominante} and the
forthcoming expressions hold for $\lambda = \pm 1$ by taking the
obvious limits. This follows from a simple continuity argument.
\end{remark}

\begin{proposition} \label{Proposicion-Leibniz}
$\Su$ and $\Di$ satisfy the following Leibniz rules
\begin{eqnarray*}
\Su (fg) \, (x) & = & \Su f(x) \, \Su g(x) + \delta(x) \, \Di f(x)
\Di g(x) \\ \Di (fg)(x) & = & \Di f(x) \Su g(x) + \Su f(x) \Di
g(x).
\end{eqnarray*}
\end{proposition}

\begin{proof}
Both expressions are straightforward and we leave them to the
reader.
\end{proof}

If $(x_t)$ is a $\Po$-sequence, then we can use our definitions of
the operators $\Su$ and $\Di$ to obtain the following relations
\begin{equation} \label{Ecuacion-Operadores}
\begin{array}{rclcl} \Su f(x_t) & = & f_s(x_{t + \frac{1}{2}},
x_{t - \frac{1}{2}}) & = & \displaystyle \frac{f(x_{t +
\frac{1}{2}}) + f(x_{t - \frac{1}{2}})}{2} \\ \Di f(x_t) & = &
f_d(x_{t + \frac{1}{2}}, x_{t - \frac{1}{2}}) & = & \displaystyle
\frac{f(x_{t + \frac{1}{2}}) - f(x_{t - \frac{1}{2}})}{x_{t +
\frac{1}{2}} - x_{t - \frac{1}{2}}}.
\end{array}
\end{equation}
Also, it is not difficult to check the validity of the following
identities
\begin{equation} \label{Ecuacion-S-D}
\begin{array}{rclcl} \Su f(x_t) & = & f(x_{t + \frac{1}{2}}) -
\displaystyle \frac{x_{t + \frac{1}{2}} - x_{t - \frac{1}{2}}}{2}
\, \Di f(x_t) \\ \Su f(x_t) & = & f(x_{t - \frac{1}{2}}) +
\displaystyle \frac{x_{t + \frac{1}{2}} - x_{t - \frac{1}{2}}}{2}
\, \Di f(x_t).
\end{array}
\end{equation}

\begin{proposition} \label{Proposicion-Iteraciones}
$\Su$ and $\Di$ satisfy the following relations
\begin{eqnarray*}
\Di \Su f(x) & = & (a+1) (A(x) - x) \, \Di^2 f(x) + a \, \Su \Di
f(x) \\ \Su^2 f(x) & = & a \, \delta(x) \, \Di^2 f(x) + (a+1)
(A(x) - x) \, \Su \Di f(x) + f(x).
\end{eqnarray*}
\end{proposition}

\begin{proof}
Given a complex number $x$, let $(x_t)$ be a $\Po$-sequence with
base point $x$. By continuity, we can assume that $$x_{t +
\frac{1}{2}} \neq x_{t - \frac{1}{2}} \qquad \mbox{for} \qquad
\textstyle t = 0, \pm \frac{1}{2}.$$ Under these assumptions, we
can use Corollary \ref{Corolario-Formulario} and
(\ref{Ecuacion-Operadores}) to obtain
\begin{eqnarray*}
\Di \Su f(x) & = & \frac{1}{x_{\frac{1}{2}} - x_{-\frac{1}{2}}}
\Big( \frac{f(x_1) + f(x)}{2} - \frac{f(x) + f(x_{-1})}{2} \Big)
\\ 4 a \, \Su \Di f(x) & = & \frac{x_1 - x_{-1}}{x_{\frac{1}{2}} -
x_{-\frac{1}{2}}} \Big( \frac{f(x_1) - f(x)}{x_1 - x} + \frac{f(x)
- f(x_{-1})}{x - x_{-1}} \Big) \\ (a+1) (A(x) - x) \, \Di^2 f(x) &
= & \frac{x_1 - 2x + x_{-1}}{4 (x_{\frac{1}{2}} -
x_{-\frac{1}{2}})} \Big( \frac{f(x_1) - f(x)}{x_1 - x} -
\frac{f(x) - f(x_{-1})}{x - x_{-1}} \Big) \\ (a+1) (A(x) - x) \,
\Su \Di f(x) & = & \frac{x_1 - 2x + x_{-1}}{8} \Big( \frac{f(x_1)
- f(x)}{x_1 - x} + \frac{f(x) - f(x_{-1})}{x - x_{-1}} \Big) \\ a
\, \delta(x) \, \Di^2 f(x) & = & \frac{x_1 - x_{-1}}{8} \Big(
\frac{f(x_1) - f(x)}{x_1 - x} - \frac{f(x) - f(x_{-1})}{x -
x_{-1}} \Big) \\ \Su^2 f(x) & = & \frac{f(x_1) + f(x)}{4} +
\frac{f(x) + f(x_{-1})}{4}.
\end{eqnarray*}
At this point it is easy to check the stated relations. This
completes the proof.
\end{proof}

\begin{remark} \label{Observacion-Classical-Forms}
Once we have introduced the operators $\Su$ and $\Di$, we are in
position to justify some terminology introduced above. It follows
from Table I that the $\Po$-sequences of the canonical form
$\mathbf{C}$ are constant. In particular, $\Su$ is the identity
operator while $\Di$ is given by the classical derivative $$\Di f
= \frac{df}{dx}.$$ This is why we have decided to call this
canonical form the continuous form. On the other hand, the
$\Po$-sequences of the arithmetic canonical form are given by $x_t
= x_0 \pm t$. In particular, $\Su$ can be regarded as an
arithmetic mean while $\Di$ can be rewritten in terms of the
classical operator $\Delta f (x) = f(x+1) - f(x)$ as follows $$\Di
f(x) = \Delta f (x - \mbox{$\frac{1}{2}$}).$$ In summary, these
canonical forms can be regarded as the \lq classical\rq${}$ forms.
\end{remark}

\begin{remark}
By Proposition \ref{Proposicion-Iteraciones} we can characterize
those polynomials $\Po$ for which we have $\Su \Di = \Di \Su$.
Namely, the operators $\Su$ and $\Di$ commute if and only if the
coefficients of $\Po$ satisfy $a=1$ and $b=0$. In other words,
this happens only for the classical forms described in Remark
\ref{Observacion-Classical-Forms}.
\end{remark}

\subsection{Taylor formula}

In order to construct the polynomials which will appear in the
Taylor formula for the Askey-Wilson operator, we shall need the
following lemma, which is a simple consequence of Proposition
\ref{Proposicion-Formulario}.

\begin{lemma} \label{Lemma-Pj-sucesiones}
Given a $\Po$-sequence $(x_t)$, let us consider the sequence $y_t
= x_{jt}$ for some integer number $j \in \Z$. Then we have
$$\frac{y_{t + \frac{1}{2}} + y_{t - \frac{1}{2}}}{2} = a_j y_t +
b_j \qquad \mbox{and} \qquad (y_{t + \frac{1}{2}} - y_{t -
\frac{1}{2}})^2 = 4 \Big( \frac{\lambda^j - \lambda^{-j}}{\lambda
- \lambda^{-1}} \Big)^2 \delta(y_t)$$ where the coefficients $a_j$
and $b_j$ are given by $$a_j = \frac{\lambda^j + \lambda^{-j}}{2}
\qquad \mbox{and} \qquad b_j = \frac{\lambda^j + \lambda^{-j} -
2}{\lambda + \lambda^{-1} - 2} \ b.$$
\end{lemma}

It turns out by Lemma \ref{Lemma-Pj-sucesiones} that $(y_t)$ is a
$\Po_j$-sequence, where $\Po_j$ is the following symmetric
polynomial $$\Po_j(x,y) = (y - a_j x - b_j)^2 - \Big(
\frac{\lambda^j - \lambda^{-j}}{\lambda - \lambda^{-1}} \Big)^2
\delta(x) \quad \mbox{with} \quad \delta_j(x) = \Big(
\frac{\lambda^j - \lambda^{-j}}{\lambda - \lambda^{-1}} \Big)^2
\delta(x).$$ The symmetric polynomials $\Po_j$ are naturally
transformed under the action of $\mbox{Aff}(\C)$. Namely, we have
$$(g \cdot \Po)_j (x,y) = (g \cdot \Po_j) (x,y).$$ Given a
positive integer $k \ge 1$, we define the \textbf{polynomials}
$\Phi_k(x,y)$ as follows
$$\Phi_{2k} (x,y) = \prod_{j=1}^k \Po_{2j-1}(x,y) \qquad
\mbox{and} \qquad \Phi_{2k+1} (x,y) = (y-x) \prod_{j=1}^k
\Po_{2j}(x,y).$$ We set $\Phi_0(x,y) = 1$ and $\Phi_1(x,y) = y-x$.
In particular, $\Phi_k(x,y) = (-1)^k \Phi_k(y,x)$ for all $k \ge
0$. Moreover, given two $\Po$-sequences $(x_t)$ and $(y_t)$, it
follows by Lemma \ref{Lemma-Pj-sucesiones} that
\begin{equation} \label{Ecuacion-Formula-Phi}
\Phi_k(x_0,y_0) = \prod_{j=0}^{k-1} (y_0 - x_{j - \frac{k-1}{2}})
= \prod_{j=0}^{k-1} (y_{j - \frac{k-1}{2}} - x_0) \qquad
\mbox{for} \quad k \ge 1.
\end{equation}

\begin{remark}
The polynomials $\Phi_k$ depend continuously on the coefficients
$a,b,c$. In particular, once again we can take limits $\lambda
\rightarrow \pm 1$ in the forthcoming expressions.
\end{remark}

\begin{lemma} \label{Lemma-DPhik}
The following relations hold for any $k \ge 1$ $$\Di \Phi_k
(\cdot,y) = \, \frac{\lambda^{-k} - \lambda^k}{\lambda -
\lambda^{-1}} \ \Phi_{k-1}(\cdot,y)$$ $$\ \Di \Phi_k (x, \cdot) =
\, \frac{\lambda^k - \lambda^{-k}}{\lambda - \lambda^{-1}} \
\Phi_{k-1} (x, \cdot).$$ Moreover, for $k = 0$ we have the
relations $\Di \Phi_0 (\cdot,y) = \Di \Phi_0 (x,\cdot) = 0$.
\end{lemma}

\begin{proof}
Given a complex number $x$, we consider a $\Po$-sequence with base
point $x$. Moreover, by continuity we can assume that $x_{1/2}
\neq x_{-1/2}$. In that case, it is not difficult to check that
$$\frac{\Phi_k(x_{1/2},y) - \Phi_k(x_{-1/2},y)}{x_{1/2} -
x_{-1/2}} = \frac{x_{-k/2} - x_{k/2}}{x_{1/2} - x_{-1/2}} \
\Phi_{k-1}(x_0,y).$$ The first relation is then obtained by
Proposition \ref{Proposicion-Formulario}. The second relation
follows by the first relation and the identity $\Phi_k(x,y) =
(-1)^k \Phi_k(y,x)$. Finally, when $k=0$ the stated relations are
obvious. This completes the proof.
\end{proof}

\begin{lemma} \label{Lema-Partial-Phi}
The linear operator $\partial_k: \Pol_d[x] \rightarrow
\Pol_{d-k}[x]$ defined by $$\partial_k f (x) = \Big(
\prod_{j=0}^{k-1} \frac{\lambda - \lambda^{-1}}{\lambda^{k-j} -
\lambda^{j-k}} \Big) \, \Di^k f (x) \qquad \mbox{for} \qquad
|\lambda| \neq 1$$ has a continuous extension for any $\lambda$
which depends polynomially on $\Po$.
\end{lemma}

\begin{proof}
From Lemma \ref{Lemma-DPhik}, the following relation holds
whenever $r \ge k$ and the parameter $\lambda$ satisfies
$|\lambda| \neq 1$
\begin{equation} \label{Ecuacion-Partial-Phi}
\partial_k \Phi_r (x,\cdot) = \Big( \prod_{j=0}^{k-1}
\frac{\lambda^{r-j} - \lambda^{j-r}}{\lambda^{k-j} -
\lambda^{j-k}} \Big) \, \Phi_{r-k} (x,\cdot).
\end{equation}
Moreover, $\partial_k \Phi_r (x,\cdot)$ vanishes for $r < k$. By a
well-known result of $q$-combinatorics, see for instance
\cite{VW}, the rational function which appears in
(\ref{Ecuacion-Partial-Phi}) is a symmetric polynomial in
$\lambda$ and $\lambda^{-1}$. In particular, the function
$\partial_k \Phi_r (x,\cdot)$ can be written as a polynomial in $a
= \frac{1}{2} (\lambda + \lambda^{-1}), b$ and $c$. Since for any
fixed $x$, the family $\Phi_0 (x,\cdot), \Phi_1 (x,\cdot), \ldots
\Phi_d (x,\cdot)$ is a basis of $\Pol_d[x]$, the proof is
completed.
\end{proof}

\begin{remark}
The operator $\partial_k$ defined in Lemma \ref{Lema-Partial-Phi}
can be regarded as the analog of the \emph{divided power of the
$k$-th derivative} $\partial^{[k]}$, which appears in finite field
theory.
\end{remark}

\begin{remark} \label{Observacion-Coeficientes-q-Binomiales}
Taking $q = \lambda^2$, the rational function appearing in the
proof of Lemma \ref{Lema-Partial-Phi} can be rewritten in terms of
a $q$-binomial coefficient $$\prod_{j=0}^{k-1} \frac{\lambda^{r-j}
- \lambda^{j-r}}{\lambda^{k-j} - \lambda^{j-k}} = \lambda^{-k
(r-k)} \Big[ \!\!
\begin{array}{c} r \\ k \end{array} \!\! \Big]_q \qquad
\mbox{where} \qquad \Big[ \!\!
\begin{array}{c} r \\ k \end{array} \!\! \Big]_q =
\frac{(q;q)_r}{(q;q)_k (q;q)_{r-k}}.$$
\end{remark}

The following result constitutes a Taylor formula for the
Askey-Wilson operator. As we shall see immediately, it generalizes
the continuous and discrete classical expressions of this formula.

\begin{theorem} \label{Teorema-Taylor-Formula}
If $(x_t)$ is a $\Po$-sequence and $f \in \Pol_r[x]$, then we have
$$f(y) = \sum_{k=0}^r \, \partial_k f(x_{k/2}) \,
\prod_{j=0}^{k-1} (y - x_j).$$
\end{theorem}

\begin{proof}
There exists a family of complex coefficients $\lambda_0,
\lambda_1, \ldots \lambda_r$ such that $$f(y) = \sum_{k=0}^r \,
\lambda_k \, \prod_{j=0}^{k-1} (y - x_j).$$ On the other hand,
since $\displaystyle \prod_{j=0}^{k-1} (y - x_j) = \Phi_k
(x_{(k-1)/2},y)$, formula (\ref{Ecuacion-Partial-Phi}) gives
$$\partial_j f(y) = \sum_{k=j}^r \lambda_k \Big( \prod_{i=0}^{j-1}
\frac{\lambda^{k-i} - \lambda^{i-k}}{\lambda^{j-i} -
\lambda^{i-j}} \Big) \, \Phi_{k-j}(x_{(k-1)/2},y).$$ By identity
(\ref{Ecuacion-Formula-Phi}), if $k > j$ then
$\Phi_{k-j}(x_{(k-1)/2},y)$ vanishes at $$y = x_{i + \frac{j}{2}}
\qquad \mbox{for} \qquad i =0,1, \ldots k-j-1.$$ In particular,
evaluating our expression for $\partial_j f(y)$ at the point $y =
x_{j/2}$, we obtain the identity $\lambda_k =
\partial_k f(x_{k/2})$ as we wanted. This completes the
proof.
\end{proof}

\begin{remark}
The classical continuous and discrete Taylor formulas are
particular cases of Theorem \ref{Teorema-Taylor-Formula}. Namely,
if $\Po(x,y) = x^2 + y^2 - 2xy$ the $\Po$-sequences are constant
and the operator $\Di$ coincides with the classical derivative, as
it was recalled in Remark \ref{Observacion-Classical-Forms}.
Therefore, Theorem \ref{Teorema-Taylor-Formula} takes the
classical form $$f(y) = \sum_{k=0}^r \frac{f^{(k)}(x)}{k!}
(y-x)^k.$$ On the other hand, if $\Po(x,y) = (x - y)^2 -
\frac{1}{4}$, we know that the $\Po$-sequences have the form $x_0
\pm t$. Thus, using the classical notation $[t]_k = t (t-1) \cdots
(t-k+1)$ and $\Delta f (x) = f(x+1) - f(x)$, we obviously have
$$\Phi_k(x,y) = [y - x + \mbox{$\frac{k-1}{2}$}]_k \qquad
\mbox{and} \qquad
\partial_k f (x) = \frac{1}{k!} \, \Delta^k f (x -
\mbox{$\frac{k}{2}$}).$$ In particular, Theorem
\ref{Teorema-Taylor-Formula} provides the discrete Taylor formula
$$f(y) = \sum_{k=0}^r \Delta^k f(x) \binom{y-x}{k}.$$
\end{remark}

\subsection{Some remarks for analytic functions}
\label{Subsection-Analytic}

Given an open subset $\Omega$ of the complex plane, we denote the
space of analytic functions in $\Omega$ by $\mathcal{H}(\Omega)$.
Also, $\gamma \simeq 0 \ (\mbox{mod} \, \Omega)$ means that
$\gamma$ is a cycle in $\Omega$ homologous to zero with respect to
$\Omega$. Finally, given $z \in \Omega$, $\mbox{Ind}(\gamma,z)$
denotes the index of $z$ with respect to $\gamma$.

\begin{lemma} \label{Lema-Cauchy}
Given an open subset $\Omega$ of the complex plane, let us denote
by $\Delta_{\Omega}$ the diagonal of $\Omega \times \Omega$. Then,
for any $f \in \mathcal{H}(\Omega)$, the function $$f_d(u,v) =
\frac{f(u)-f(v)}{u-v} \qquad \mbox{for} \qquad (u,v) \in (\Omega
\times \Omega) \setminus \Delta_{\Omega}$$ can be continuously
extended to an analytic function $f_d: \Omega \times \Omega
\rightarrow \C$. Moreover, if $\gamma \simeq 0 \ (\textnormal{mod}
\, \Omega)$ and $\textnormal{Ind}(\gamma,u) =
\textnormal{Ind}(\gamma,v) = 1$, then we have $$f_d(u,v) =
\frac{1}{2 \pi i} \int_{\gamma} \frac{f(y)}{(y-u)(y-v)} \, dy.$$
\end{lemma}

\begin{proof}
It is a simple consequence of Cauchy's integral formula.
\end{proof}

Let us consider the set $\Omega^{(1)} = \{x \in \C: A(x) \pm
\sqrt{\delta(x)} \in \Omega\}$. By Lemma \ref{Lema-Cauchy}, we can
define the operator $\Di: \mathcal{H}(\Omega) \rightarrow
\mathcal{H}(\Omega^{(1)})$ as follows
\begin{equation} \label{Ecuacion-D-Analitico}
\Di f(x) = f_d \big( A(x) + \sqrt{\delta(x)},A(x) -
\sqrt{\delta(x)} \big) = \frac{1}{2 \pi i} \int_{\gamma}
\frac{f(y)}{\Po(x,y)} \, dy,
\end{equation}
with $\gamma \simeq 0 \ (\mbox{mod} \, \Omega)$ and
$\mbox{Ind}(\gamma, A(x) \pm \sqrt{\delta(x)}) = 1$. It is clear
that (\ref{Ecuacion-D-Analitico}) extends the original definition
of the Askey-Wilson operator $\Di$. Moreover, we can also extend
the definition of the companion operator $\Su$. Namely, if $f \in
\mathcal{H}(\Omega)$, $x_0 \in \Omega^{(1)}$ and $\delta(x_0) \neq
0$, the function $$\Su f (x) = f_s \big( A(x) + \sqrt{\delta(x)},
A(x) - \sqrt{\delta(x)} \big) = \frac{f \big( A(x) +
\sqrt{\delta(x)} \big) + f \big( A(x) - \sqrt{\delta(x)}
\big)}{2}$$ is obviously analytic for $x$ in a neighborhood of
$x_0$. Besides, $\Su f$ is continuous at $x_0$ when $\delta(x_0) =
0$. Thus, the roots of $\delta$ are removable singularities and
$\Su f$ becomes analytic in $\Omega^{(1)}$. If we define
recursively the sets $$\Omega^{(k+1)} = \{x \in \C: A(x) \pm
\sqrt{\delta(x)} \in \Omega^{(k)}\} \qquad \mbox{with} \qquad
\Omega^{(0)} = \Omega,$$ we can consider the iterated operators
$\Su^k, \Di^k: \mathcal{H}(\Omega) \rightarrow
\mathcal{H}(\Omega^{(k)})$. An open subset $\Omega$ of the complex
plane will be called \textbf{$\Po$-invariant} if, for any
$\Po$-sequence with base point $x_0 \in \Omega$, we have $x_{\pm
1/2} \in \Omega$. In other words, if $\Omega \subset \Omega^{(k)}$
for any positive integer $k \ge 1$. If $\Omega$ is
$\Po$-invariant, we deduce that $\Su^k f \in \mathcal{H}(\Omega)$
and $\Di^k f \in \mathcal{H}(\Omega)$ for any function $f$
analytic in $\Omega$ and any positive integer $k$.

\begin{remark}
Similarly, if $\mathcal{M}(\Omega)$ stands for the space of
meromorphic functions in $\Omega$, we can extend the previous
operators so that $$\Su^k, \Di^k: \mathcal{M}(\Omega) \rightarrow
\mathcal{M}(\Omega^{(k)}).$$ Furthermore, if $\Omega$ is
$\Po$-invariant, then we have operators $\Su^k, \Di^k:
\mathcal{M}(\Omega) \rightarrow \mathcal{M}(\Omega)$.
\end{remark}

\begin{remark}
Obviously, Propositions \ref{Proposicion-Leibniz} and
\ref{Proposicion-Iteraciones} remain valid in this new context.
\end{remark}

\begin{proposition} \label{Proposicion-Dk-Analitico}
Let $\Omega$ be an open set of the complex plane and $\gamma
\simeq 0 \ (\textnormal{mod} \, \Omega)$. Then, given $f \in
\mathcal{H}(\Omega)$, an integer $k \ge 0$ and $x \in
\Omega^{(k)}$, we have $$\Di^k f (x) = \Big( \prod_{j=0}^{k-1}
\frac{\lambda^{k-j} - \lambda^{j-k}}{\lambda - \lambda^{-1}} \Big)
\, \frac{1}{2 \pi i} \int_{\gamma} \frac{f(y)}{\Phi_{k+1}(x,y)} \,
dy$$ if $x$ is the base point of a $\Po$-sequence $(x_t)$ with
$\textnormal{Ind}(\gamma, x_{j - \frac{k}{2}}) = 1$ for $j = 0,1,
\ldots k$.
\end{proposition}

\begin{proof}
The cases $k=0$ and $k=1$ are Cauchy's integral formula and
(\ref{Ecuacion-D-Analitico}) respectively. The general case
follows by induction from (\ref{Ecuacion-D-Analitico}) and the
relation $$\Di \Big[ \frac{1}{\Phi_k(\cdot,y)} \Big] (x) =
\frac{\lambda^k - \lambda^{-k}}{\lambda - \lambda^{-1}}
\Phi_{k+1}(x,y)^{-1}$$ which can be checked by the reader. This
completes the proof.
\end{proof}

\begin{corollary} \label{Corolario-Partialk-Analitico}
Let $\Omega$ be an open subset of the complex plane and $\gamma
\simeq 0 \ (\textnormal{mod} \, \Omega)$. Then, given $f \in
\mathcal{H}(\Omega)$, the expression $$\partial_k f(x) =
\frac{1}{2 \pi i} \int_{\gamma} \frac{f(y)}{\Phi_{k+1}(x,y)} \,
dy$$ defines a linear operator $\partial_k: \mathcal{H}(\Omega)
\rightarrow \mathcal{H}(\Omega^{(k)})$. Moreover, we have
$$\partial_k f(x) = \sum_{j=0}^k \textnormal{Res}_{\mid y = x_{j -
\frac{k}{2}}} \Big( \frac{f(y)}{\Phi_{k+1}(x,y)} \Big).$$
\end{corollary}

\begin{remark}
It turns out that Lemma \ref{Lema-Partial-Phi} is now a
consequence of Corollary \ref{Corolario-Partialk-Analitico}.
\end{remark}

\begin{remark} \label{Observacion-Display1}
When $x$ is the base point of a $\Po$-sequence $(x_t)$ with $x_{j
- \frac{k}{2}}$ pairwise distinct for $j = 0,1, \ldots k$, we have
by Corollary \ref{Corolario-Partialk-Analitico} $$\partial_k f(x)
= \sum_{j=0}^k \frac{f(x_{j - \frac{k}{2}})}{\prod_{0 \le i \neq j
\le k} (x_{j - \frac{k}{2}} - x_{i - \frac{k}{2}})}.$$ As we
explain in Remark \ref{Observacion-Display2}, this is useful to
display our Rodrigues formulas given by Theorem
\ref{Teorema-Rodrigues}. We recall that Nikiforov, Suslov and
Uvarov give similar expressions for less intrinsic operators in
\cite[3.2.3]{NSU}.
\end{remark}

\begin{remark}
Under the hypothesis of Remark \ref{Observacion-Display1}, Theorem
\ref{Teorema-Taylor-Formula} gives
$$f(x) = \sum_{k=0}^r \Big( \sum_{j=0}^k \frac{f(x_j)}{\prod_{i
\neq j} (x_j - x_i)} \Big) \prod_{j=0}^{k-1} (x - x_j)$$ which is
Newton's divided difference formula for the interpolation
polynomial. In fact, as the referee of this paper has pointed out,
the formula given in Remark \ref{Observacion-Display1} can be
proved by induction. In particular, it is possible to obtain an
alternative proof of Theorem \ref{Teorema-Taylor-Formula} from the
theory of interpolation polynomials. We prefer our proof since it
shows how the polynomials $\Phi_k$ are related to $\Di$, see Lemma
\ref{Lemma-DPhik}.
\end{remark}

\begin{remark}
The operators $\partial_k: \mathcal{H}(\Omega) \rightarrow
\mathcal{H}(\Omega^{(k)})$ provide the natural framework to
develop a generalization of Theorem \ref{Teorema-Taylor-Formula}.
Namely, in a separate work \cite{MP2}, we study the convergence of
the Taylor series associated to the operators $\partial_k$ of a
given function $f \in \mathcal{H}(\Omega)$ and its relation with
the basic hypergeometric functions.
\end{remark}

\section{Taylor coefficients and Rodrigues formula}
\label{Section3}

In this section we introduce the corresponding hypergeometric
operator. Then we compute the Taylor coefficients of its
eigenfunctions. Also, we provide a general procedure to obtain
Rodrigues type formulas. At the end of this section, we shall
study the discrete orthogonal relations which arise from our
techniques.

\subsection{The hypergeometric operator}

Given a symmetric polynomial $\Po$ in the set $\mathcal{P}$, we
consider the associated Askey-Wilson operator $\Di$ and its
companion $\Su$. Let $\sigma \in \Pol_2[x]$ and $\tau \in
\Pol_1[x]$ given by $$\sigma(x) = \alpha_2 x^2 + \alpha_1 x +
\alpha_0 \qquad \mbox{and} \qquad \tau(x) = \beta_1 x + \beta_0.$$
The aim of this section is to study the polynomic eigenfunctions
of the following operator $$\OH = \sigma \Di^2 + \tau \Su \Di.$$
$\OH$ will be called the \textbf{hypergeometric operator}
associated to $\Po$, $\sigma$ and $\tau$. In Section
\ref{Section4} we shall study in detail some relevant particular
cases. It follows by Proposition
\ref{Proposicion-TerminoDominante} that $\OH(\Pol_k[x]) \subset
\Pol_k[x]$ for $k \ge 1$. Moreover, if $f(x) = x^k$ we obtain
$$\OH f (x) = - \mu_k x^k + \ldots$$ where the dots stand for
terms of lower degree and $$\mu_k = - \frac{\lambda^k -
\lambda^{-k}}{(\lambda - \lambda^{-1})^2} \Big[ \Big( \alpha_2 +
\frac{\lambda - \lambda^{-1}}{2} \beta_1 \Big) \lambda^{k-1} -
\Big( \alpha_2 - \frac{\lambda - \lambda^{-1}}{2} \beta_1 \Big)
\lambda^{1-k} \Big].$$ As an endomorphism of $\Pol_k[x]$, the
hypergeometric operator $\OH$ depends continuously on $\Po$,
$\sigma$ and $\tau$. Therefore, taking suitable limits for
$\lambda = \pm 1$, it turns out that $\mu_k$ is a continuous
function on $\lambda$. The following formula, for which we should
take limits in $\lambda$ when $\lambda = \pm 1$, will be useful in
Section \ref{Section4} $$\mu_k - \mu_j = \frac{\lambda^{j-k} -
\lambda^{k-j}}{(\lambda - \lambda^{-1})^2} \Big[ \Big( \alpha_2 +
\frac{\lambda - \lambda^{-1}}{2} \beta_1 \Big) \lambda^{j+k-1} -
\Big( \alpha_2 - \frac{\lambda - \lambda^{-1}}{2} \beta_1 \Big)
\lambda^{1-j-k} \Big].$$ Now we point out some relevant remarks on
the operator $\OH$:
\begin{itemize}
\item If $\OH^{\Po}$ stands for the hypergeometric operator
associated to the symmetric polynomial $\Po$ and we consider some
$g = (\eta, \zeta) \in \mathfrak{T}(\C) \rtimes GL(\C) =
\mbox{Aff}(\C)$, it follows easily from Remark
\ref{Observacion-SD-Accion} that $$g \ \OH^{\Po} g^{-1} = \zeta^2
\sigma^g (\Di^{g \cdot \Po})^2 + \zeta \tau^g \Su^{g \cdot \Po}
\Di^{g \cdot \Po}$$ where $(g \cdot f)(x) = f (g^{-1} (x))$,
$\sigma^g(x) = \sigma(g^{-1}(x))$ and $\tau^g(x) =
\tau(g^{-1}(x))$.

\item When $(x_t)$ is a $\Po$-sequence with $x_{t - \frac{1}{2}}
\neq x_{t + \frac{1}{2}}$ for $t = 0, \pm \frac{1}{2}$, we can
write $$\begin{array}{r} \displaystyle \OH f (x_0) =
\frac{\sigma(x_0)}{x_{1/2} - x_{-1/2}} \Big[ \frac{f(x_1) -
f(x_0)}{x_1 - x_0} - \frac{f(x_0) - f(x_{-1})}{x_0 - x_{-1}} \Big]
\ \\ \displaystyle + \frac{\tau(x_0)}{2} \Big[ \frac{f(x_1) -
f(x_0)}{x_1 - x_0} + \frac{f(x_0) - f(x_{-1})}{x_0 - x_{-1}}
\Big]. \end{array}$$ In fact, there exist a limiting version of
this equality for any $\Po$-sequence $(x_t)$. Moreover, since $x_1
- x_{-1} = 2 a (x_{1/2} - x_{-1/2})$ by Corollary
\ref{Corolario-Formulario}, $\OH f (x_0)$ can be written in terms
of $x_t$ with $t \in \Z$. In other words, since the subsequence of
$(x_t)$ which arise when $t$ runs over $\Z$ is a $\Po_2$-sequence
with the terminology introduced after Lemma
\ref{Lemma-Pj-sucesiones}, we can say that $\OH$ only depends on
$\Po_2$.

\item If the symmetric polynomial $\Po$ is of type \textbf{O}
(resp. \textbf{E}), it can be checked that $\Po_2$ is of type
\textbf{A} (resp. \textbf{C}). Besides, if $\Po$ is of type
\textbf{A} (resp. \textbf{C}), we also have that $\Po_2$ is of
type \textbf{A} (resp. \textbf{C}). Therefore, taking into account
the fact that $\OH$ only depends on $\Po_2$, we do not loose
generality if we restrict our study to the canonical forms
\textbf{T}, \textbf{G}, \textbf{Q}, \textbf{A} and \textbf{C}.
This point of view leads to a different and somehow more natural
classification than the one given in \cite[3.4]{NSU}. Namely,
while our classification is purely complex, it can be said that
the classification given in \cite{NSU} lives over the real field.
As we point out below, another remarkable advantage of our
formulation is that affine transformations become trivial.

\item If $\lambda = \pm i$, we have $\Di^2 = (\lambda +
\lambda^{-1})\partial_2 = 0$. In particular, $\OH = \tau \Su \Di$
and so we get a degenerate case. This phenomenon can be avoided by
considering the more general operator $$\OH' = \sigma
\partial_2 + \tau \Su \Di.$$ It is not difficult to adapt the
methods we shall develop to this case. When $\lambda = \pm i$ in
$\Po$, it can be checked that $\Po_2$ is of type \textbf{E}.
Therefore, the even canonical form can be regarded from this
alternative point of view.

\item The polynomial $\Po_2$ can not be of type
\textbf{O}. However, we can construct a limiting operator which
somehow corresponds to certain $\Po_2$ of type \textbf{O}. Namely,
it suffices to take a sequence $\Po(n)$ in $\mathcal{P}$ such that
the sequence $\Po(n)_2$ converges in $\mathcal{P}$ to a symmetric
polynomial of type \textbf{O}. It turns out that, in this
situation, the sequence $\Po(n)$ does not converge when $n
\rightarrow \infty$. However, it can be checked that the
corresponding hypergeometric operators $\OH(n)$ do converge.
Therefore, the oscillating canonical form could be analyzed by
means of the operator $$\OH = \lim_{n \rightarrow \infty}
\OH(n).$$ Although it can be studied what happens with our methods
after taking this limit, we shall not cover this problem here.
\end{itemize}

\subsection{Taylor coefficients and Rodrigues formula}

The first step in our process is the following commutation
relation for the iterated hypergeometric operators. Expressions of
the same kind can also be found in \cite{Ba}.

\begin{lemma} \label{Lema-Conmutacion}
Taking $\sigma_0 = \sigma$, $\tau_0 = \tau$ and given $j \ge 0$,
let us define the polynomials $\sigma_j \in \Pol_2[x]$ and $\tau_j
\in \Pol_1[x]$ inductively as follows
\begin{eqnarray*}
\sigma_{j+1}(x) & = & \Su \sigma_j (x) + (a+1) (A(x) - x) \Su
\tau_j (x) + a \delta(x) \Di \tau_j(x) \\ \tau_{j+1}(x) & = & \Di
\sigma_j (x) + a \Su \tau_j(x) + (a+1) (A(x) - x) \Di \tau_j(x)
\end{eqnarray*}
Then, the iterated hypergeometric operator $\OH_j = \sigma_j \Di^2
+ \tau_j \Su \Di$ satisfies $$\Di (\OH_j - \mu_j) = (\OH_{j+1} -
\mu_{j+1}) \Di.$$
\end{lemma}

\begin{proof}
Applying the Leibniz rules and the commuting relations developed
for the operators $\Su$ and $\Di$ in Propositions
\ref{Proposicion-Leibniz} and \ref{Proposicion-Iteraciones}, we
easily obtain
\begin{eqnarray*}
\Di \OH_j f (x) & \!\! = \!\! & \Di (\sigma_j \Di^2 f) (x) + \Di
(\tau_j \Su \Di f) (x) \\ & \!\! = \!\! & \Su \sigma_j(x) \Di^3 f
(x) + \Di \sigma_j(x) \Su \Di^2 f (x) + \Su \tau_j (x) \Di \Su \Di
f (x) + \Di \tau_j (x) \Su^2 \Di f (x) \\ & \!\! = \!\! & \Su
\sigma_j(x) \Di^3 f (x) + \Di \sigma_j(x) \Su \Di^2 f (x) \\ &
\!\! + \!\! & \Su \tau_j (x) \big( (a+1) (A(x)-x) \Di^3 f (x) + a
\Su \Di^2 f (x) \big) \\ & \!\! + \!\! & \Di \tau_j (x) \big( a
\delta(x) \Di^3 f (x) + (a+1) (A(x)-x) \Su \Di^2 f (x) + \Di f (x)
\big) \\ & \!\! = \!\! & \sigma_{j+1}(x) \Di^3 f (x) +
\tau_{j+1}(x) \Su \Di^2 f (x) + \Di \tau_j (x) \Di f (x).
\end{eqnarray*}
In other words, $\Di \OH_j f (x) = \big( \OH_{j+1} + \Di \tau_j
\big) \Di f (x)$ where $\Di \tau_j$ is a constant $\nu_j$ since
$\tau_j \in \Pol_1[x]$. On the other hand, taking $f(x) = x^k$, we
define $\mu_{j,k}$ by the relation $$\OH_j f (x) = - \mu_{j,k} x^k
+ \cdots$$ where the dots stand for terms of lower degree. Then,
Proposition \ref{Proposicion-TerminoDominante} gives $$- \Big(
\frac{\lambda^k - \lambda^{-k}}{\lambda - \lambda^{-1}} \Big)
\mu_{j,k} = \Big( \frac{\lambda^k - \lambda^{-k}}{\lambda -
\lambda^{-1}} \Big) (- \mu_{j+1,k-1} + \nu_j).$$ By continuity, we
have $\mu_{j+1,k-1} = \mu_{j,k} + \nu_j$ or equivalently
$$\mu_{j,k} = \mu_{j+k} + \sum_{i=0}^{j-1} \nu_i.$$ In particular,
we obtain that $\mu_j + \nu_0 + \ldots + \nu_{j-1} = \mu_{j,0} =
0$. Therefore, we deduce the identity $\nu_j = \mu_j - \mu_{j+1}$
as we wanted. This completes the proof.
\end{proof}

Given a symmetric polynomial $\Po(x,y) = x^2 + y^2 -2axy - 2b(x+y)
+ c$, let us take a $\Po$-function $\theta_t$. Then we define the
\textbf{auxiliary functions} $$\widetilde{\sigma}_j^{\pm} (t) =
\sigma_j (\theta_t) \pm \frac{\theta_{t + \frac{1}{2}} - \theta_{t
- \frac{1}{2}}}{2} \tau_j(\theta_t) \qquad \mbox{for} \qquad j \ge
0.$$ We shall write $\widetilde{\sigma}^{\pm}$ to denote the functions
$\widetilde{\sigma}_0^{\pm}$ respectively.

\begin{lemma} \label{Lema-Funciones-Auxiliares}
The auxiliary functions satisfy the following relations
$$\textstyle \widetilde{\sigma}_j^+ (t) = \widetilde{\sigma}^+ (t
+ \frac{j}{2}) \qquad \mbox{and} \qquad \widetilde{\sigma}_j^- (t)
= \widetilde{\sigma}^- (t - \frac{j}{2}).$$
\end{lemma}

\begin{proof}
By definition we have
\begin{eqnarray*}
\widetilde{\sigma}_{j+1}^{\pm} (t) & = & \sigma_{j+1}(\theta_t) \pm
\textstyle \frac{1}{2} (\theta_{t + \frac{1}{2}} - \theta_{t -
\frac{1}{2}}) \displaystyle \tau_{j+1}(\theta_t) \\ & = & \Su
\sigma_j (\theta_t) + (a+1)(A(\theta_t) - \theta_t) \Su
\tau_j(\theta_t) + a \delta(\theta_t) \Di \tau_j (\theta_t) \\ &
\pm & \textstyle \frac{1}{2} (\theta_{t + \frac{1}{2}} - \theta_{t
- \frac{1}{2}}) \displaystyle \big( \Di \sigma_j(\theta_t) + a \Su
\tau_j(\theta_t) + (a+1) (A(\theta_t) - \theta_t) \Di \tau_j
(\theta_t) \big).
\end{eqnarray*}
Then, by Corollary \ref{Corolario-Formulario} and
(\ref{Ecuacion-Operadores}), we get
\begin{eqnarray*}
\widetilde{\sigma}_{j+1}^{\pm} (t) & = & \frac{\sigma_j(\theta_{t +
\frac{1}{2}}) + \sigma_j(\theta_{t - \frac{1}{2}})}{2} +
\frac{\theta_{t+1} - 2 \theta_t + \theta_{t-1}}{8} \big(
\tau_j(\theta_{t + \frac{1}{2}}) + \tau_j(\theta_{t -
\frac{1}{2}}) \big) \\ & + & \frac{\theta_{t+1} - \theta_{t-1}}{8}
\big( \tau_j(\theta_{t + \frac{1}{2}}) - \tau_j(\theta_{t -
\frac{1}{2}}) \big) \pm \frac{\theta_{t+1} - \theta_{t-1}}{8}
\big( \tau_j(\theta_{t + \frac{1}{2}}) + \tau_j(\theta_{t -
\frac{1}{2}}) \big) \\ & \pm & \frac{\theta_{t+1} - 2 \theta_t +
\theta_{t-1}}{8} \big( \tau_j(\theta_{t + \frac{1}{2}}) -
\tau_j(\theta_{t - \frac{1}{2}}) \big) \pm
\frac{\sigma_j(\theta_{t + \frac{1}{2}}) - \sigma_j(\theta_{t -
\frac{1}{2}})}{2}.
\end{eqnarray*}
Simplifying the expression above, we obtain
$\widetilde{\sigma}_{j+1}^{\pm} (t) = \widetilde{\sigma}_j^{\pm}(t \pm
\frac{1}{2})$. Therefore, the stated relations follow easily by
induction and so the proof is completed.
\end{proof}

Now, given a positive integer $d$, we assume that $\mu_j \neq
\mu_k$ for $0 \le j < k \le d$. In particular, this condition
implies that $\OH$ is a diagonalizable operator in $\Pol_d[x]$.
Therefore, for any $0 \le k \le d$, there exists an eigenfunction
$f_k \in \Pol_k[x]$ with degree $k$ satisfying $$\OH f_k (x) +
\mu_k f_k (x) = 0.$$ Moreover, if we define $f_k^{(j)} = \Di^j
f_k$, then Lemma \ref{Lema-Conmutacion} and induction gives
\begin{equation} \label{Ecuacion-Autofunciones-Iteradas}
\OH_j f_k^{(j)} + (\mu_k - \mu_j) f_k^{(j)} = 0.
\end{equation}

\begin{lemma} \label{Lema-Ecuacion-Recurrente}
The following relation holds for any $j \ge 0$
$$\widetilde{\sigma}_j^-(t) f_k^{(j+2)} (\theta_t) +
\tau_j(\theta_t) f_k^{(j+1)} (\theta_{t + \frac{1}{2}}) + (\mu_k -
\mu_j) f_k^{(j)} (\theta_t) = 0.$$
\end{lemma}

\begin{proof}
Applying the definition of $\widetilde{\sigma}_j^-$ in terms of
$\sigma_j$ and $\tau_j$ and relation (\ref{Ecuacion-S-D}), we can
rewrite the left hand side as $$\mathrm{LHS} = \sigma_j(\theta_t)
\Di^2 f_k^{(j)} (\theta_t) + \tau_j (\theta_t) \Su \Di f_k^{(j)}
(\theta_t) + (\mu_k - \mu_j) f_k^{(j)} (\theta_t).$$ In other
words, the left hand side coincides with $\OH_j
f_k^{(j)}(\theta_t) + (\mu_k - \mu_j) f_k^{(j)}(\theta_t)$, which
vanishes by relation (\ref{Ecuacion-Autofunciones-Iteradas}).
Therefore, the proof is completed.
\end{proof}

If we take the $\Po$-function $\theta_t$ in such a way that
$\widetilde{\sigma}^-(0) = 0$, then we can give explicit
expressions for the eigenfunctions $f_0, f_1, \ldots f_d$ in terms
of the Taylor formula developed in Section \ref{Section2}.

\begin{theorem} \label{Teorema-Expresiones-Polinomicas}
Let us assume that $\widetilde{\sigma}^-(0) = 0$, then the
eigenfunctions $f_0, f_1, \ldots f_d$ of $\OH$ satisfy the
following identities $$f_k (x) =
\partial_k f_k \sum_{j=0}^k \Big( \prod_{i=j}^{k-1}
\frac{\lambda^{1+i} - \lambda^{-1-i}}{\lambda - \lambda^{-1}} \,
\frac{\tau_i(\theta_{i/2})}{\mu_i - \mu_k} \Big) \prod_{i=0}^{j-1}
(x - \theta_i)$$ where $\partial_k f_k$ stands for the coefficient
of $x^k$ in $f_k(x)$. Moreover, assuming $f_k(\theta_0) \neq 0$,
we also have the following identities $$f_k (x) = f_k(\theta_0)
\sum_{j=0}^k \Big( \prod_{i=0}^{j-1} \frac{\lambda -
\lambda^{-1}}{\lambda^{1+i} - \lambda^{-1-i}} \, \frac{\mu_i -
\mu_k}{\tau_i(\theta_{i/2})} \Big) \prod_{i=0}^{j-1} (x -
\theta_i).$$
\end{theorem}

\begin{proof}
Since $\widetilde{\sigma}^-(0) = 0$, we obtain from Lemmas
\ref{Lema-Funciones-Auxiliares} and \ref{Lema-Ecuacion-Recurrente}
the identity $$f_k^{(j)} (\theta_{j/2}) =
\frac{\tau_j(\theta_{j/2})}{\mu_j - \mu_k} f_{k}^{(j+1)}
(\theta_{(j+1)/2}).$$ If we write this recurrence in terms of the
divided derivative operator, we get $$\partial_j f_k
(\theta_{j/2}) = \frac{\lambda^{1+j} - \lambda^{-1-j}}{\lambda -
\lambda^{-1}} \, \frac{\tau_j(\theta_{j/2})}{\mu_j - \mu_k} \,
\partial_{j+1} f_{k} (\theta_{(j+1)/2}).$$ Then, the first relation
follows by iterating this recurrence and applying Theorem
\ref{Teorema-Taylor-Formula}. The second relation is an obvious
consequence of the first relation and the identity $$f_k
(\theta_0) = \partial_k f_k \prod_{i=0}^{k-1} \frac{\lambda^{1+i}
- \lambda^{-1-i}}{\lambda - \lambda^{-1}} \,
\frac{\tau_i(\theta_{i/2})}{\mu_i - \mu_k},$$ which follows from
the expression above for $j = 0$. This completes the proof.
\end{proof}

\begin{remark}
Although using a less clear notation, Nikiforov, Suslov and Uvarov
gave in \cite[3.1.27]{NSU} an equivalent identity to Lemma
\ref{Lema-Funciones-Auxiliares}. Moreover, Lemma
\ref{Lema-Ecuacion-Recurrente} is quite close to
\cite[3.1.24]{NSU}. However, since there is no Taylor formula for
the Askey-Wilson operators in \cite{NSU}, it seems that the
authors did not see the relevance of these results to obtain
explicit polynomic expressions for the eigenfunctions.
\end{remark}

\begin{remark}
Replacing $\theta_t$ by the $\Po$-function $\theta'_t =
\theta_{\pm(t - t_0)}$, we can use Theorem
\ref{Teorema-Expresiones-Polinomicas} to obtain explicit
expressions of the eigenfunctions $f_k$ for each $t_0$ satisfying
either $\widetilde{\sigma}^+(t_0) = 0$ or
$\widetilde{\sigma}^-(t_0) = 0$. Such points $t_0$ are those for
which $\theta_{t_0}$ is a root of the polynomial $\mathrm{Q}(x) =
\sigma(x)^2 - \delta(x) \tau(x)^2$. Since the degree of
$\mathrm{Q}$ is less or equal than 4, we obtain in the generic
cases four explicit expressions. We shall see in Section
\ref{Section4} that there are only two cases, with distinct
eigenvalues, for which $\mathrm{Q}$ has no roots.
\end{remark}

\begin{remark} \label{Observacion-Diagonalizabilidad}
The assumption $\mu_j \neq \mu_k$ for $0 \le j < k \le d$ on the
eigenvalues is too restrictive in Theorem
\ref{Teorema-Expresiones-Polinomicas}. Namely, assuming that
$\widetilde{\sigma}^-(0) = 0$, it can be proved that $\OH$ is
diagonalizable in $\Pol_d[x]$ if and only if $$\prod_{i=j}^{k-1}
\frac{\lambda^{1+i} - \lambda^{-1-i}}{\lambda - \lambda^{-1}}
\tau_i(\theta_{i/2}) = 0$$ whenever $\mu_j = \mu_k$ for some $0
\le j < k \le d$. In this case, the first identity of Theorem
\ref{Teorema-Expresiones-Polinomicas} holds summing from $j_k = 1
+ \max \{j: \, \mu_j = \mu_k\}$. In any case, for the sake of
clarity, we shall assume that the eigenvalues are pairwise
distinct.
\end{remark}

The first relation in Theorem
\ref{Teorema-Expresiones-Polinomicas} is the appropriate one when
we normalize the eigenfunctions $f_k$ so that they become monic
polynomials. The second relation in Theorem
\ref{Teorema-Expresiones-Polinomicas} will be very useful in
Section \ref{Section4}, since it will allow us to express the
eigenfunctions $f_k$ in terms of basic hypergeometric functions.
The next result can be regarded as an algorithm to obtain
Rodrigues type formulas in this setting. We recall that the
$\Po$-invariant domains were already introduced in Paragraph
\ref{Subsection-Analytic}.

\begin{lemma} \label{Lema-Equivalencias}
Let $\rho$ and $\rho_1$ be meromorphic functions in $\Omega$ for
some $\Po$-invariant open set $\Omega$. Then, the following are
equivalent:
\begin{itemize}
\item[\textnormal(a)] For any function $f \in
\mathcal{M}(\Omega)$, we have $$\rho(x) \OH f (x) = \Di (\rho_1
\Di f) (x).$$
\item[\textnormal(b)] The functions $\rho$ and $\rho_1$ satisfy
the relations $$\rho(x) \sigma(x) = \Su \rho_1 (x) \qquad
\mbox{and} \qquad \rho(x) \tau(x) = \Di \rho_1 (x).$$
\item[\textnormal(c)] The functions $\rho$ and $\rho_1$ satisfy
the relations $$\begin{array}{rl} \mathrm{(c1)} & \Di (\rho
\sigma) (x) = (a+1)(A(x)-x) \Di (\rho \tau) (x) + a \Su (\rho
\tau) (x). \\ \mathrm{(c2)} & \rho_1 (x) = \Su (\rho \sigma) (x) -
a \delta(x) \Di (\rho \tau) (x) - (a+1)(A(x)-x) \Su (\rho \tau)
(x). \ \end{array}$$
\end{itemize}
Moreover, if the $\Po$-function $\theta_t$ is not $1$-periodic,
they are also equivalent to
\begin{itemize}
\item[\textnormal(d)] The functions $\rho$ and $\rho_1$ satisfy
the relations $$\begin{array}{rl} \mathrm{(d1)} & \rho (\theta_{t
- \frac{1}{2}}) \widetilde{\sigma}^+ (t - \frac{1}{2}) =
\rho(\theta_{t + \frac{1}{2}}) \widetilde{\sigma}^-(t +
\frac{1}{2}). \\ \mathrm{(d2)} & \rho_1 (\theta_t) = \textstyle
\frac{1}{2} \big( \rho(\theta_{t - \frac{1}{2}})
\widetilde{\sigma}^+(t - \frac{1}{2}) + \rho(\theta_{t +
\frac{1}{2}}) \widetilde{\sigma}^-(t + \frac{1}{2}) \big). \qquad
\qquad \qquad
\end{array}$$
\end{itemize}
Besides, the functional equation \textnormal{(c1)} $($resp.
\textnormal{(c2)}$)$ is equivalent to \textnormal{(d1)} $($resp.
\textnormal{(d2)}$)$. Finally, if \textnormal{(a)} holds then
$\rho_1$ satisfies \textnormal{(c1)} and \textnormal{(d1)} for the
iterated operator $\OH_1$.
\end{lemma}

\begin{proof}
Properties (a) and (b) are equivalent since, by Proposition
\ref{Proposicion-Leibniz}, we have
\begin{eqnarray*}
\rho(x) \OH f (x) & = & \rho(x) \sigma(x) \Di^2 f(x) + \rho(x)
\tau(x) \Su \Di f(x) \\ \Di (\rho_1 \Di f) (x) & = & (\Su \rho_1)
(x) \, \Di^2 f(x) + (\Di \rho_1) (x) \Su \Di f(x).
\end{eqnarray*}
To see $(\textnormal{b} \Rightarrow \textnormal{c})$, we apply
Proposition \ref{Proposicion-Iteraciones}
\begin{eqnarray*}
\Di (\rho \sigma) (x) & = & (a+1)(A(x)-x) \Di^2 \rho_1 + a \Su \Di
\rho_1 (x) \\ \Su (\rho \sigma) (x) & = & a \delta(x) \Di(\rho
\tau) (x) + (a+1)(A(x)-x) \Su (\rho \tau) (x) + \rho_1 (x).
\end{eqnarray*}
Let us denote by L(c1) (resp. R(c1)) the left (resp. right) hand
side of (c1). For the equivalence between (c1) and (d1), we use
Corollary \ref{Corolario-Formulario} and
(\ref{Ecuacion-Operadores}). Namely, arguing as in Proposition
\ref{Proposicion-Iteraciones}, we easily get $$\big(
\mathrm{L(c1)} - \mathrm{R(c1)} \big) (\theta_t) =
\frac{\widetilde{\sigma}^-(t + \frac{1}{2}) \rho(\theta_{t +
\frac{1}{2}}) - \widetilde{\sigma}^+(t - \frac{1}{2})
\rho(\theta_{t - \frac{1}{2}})}{\theta_{t + \frac{1}{2}} -
\theta_{t -\frac{1}{2}}}.$$ The equivalence between (c2) and (d2)
follows by similar arguments. Now, let us assume that property (d)
holds. Then we will have
\begin{eqnarray*}
\rho_1(\theta_{t + \frac{1}{2}}) & = & \rho(\theta_t)
\widetilde{\sigma}^+(t) = \rho(\theta_t) \Big( \sigma(\theta_t) +
\frac{\theta_{t + \frac{1}{2}} - \theta_{t - \frac{1}{2}}}{2}
\tau(\theta_t) \Big) \\ \rho_1(\theta_{t - \frac{1}{2}}) & = &
\rho(\theta_t) \widetilde{\sigma}^-(t) = \rho(\theta_t) \Big(
\sigma(\theta_t) - \frac{\theta_{t + \frac{1}{2}} - \theta_{t -
\frac{1}{2}}}{2} \tau(\theta_t) \Big).
\end{eqnarray*}
In particular, $\Su \rho_1 (\theta_t) = \rho(\theta_t)
\sigma(\theta_t)$ and $\Di \rho_1 (\theta_t) = \rho(\theta_t)
\tau(\theta_t)$. That is, we have seen $(\textnormal{d}
\Rightarrow \textnormal{b})$. It remains to see that (a) implies
that $\rho_1$ satisfies (d1) for $\OH_1$. Since (a) implies (d),
we have by Lemma \ref{Lema-Funciones-Auxiliares} $$\textstyle
\rho_1(\theta_{t - \frac{1}{2}}) \widetilde{\sigma}_1^+(t -
\frac{1}{2}) = \rho(\theta_t) \widetilde{\sigma}^-(t)
\widetilde{\sigma}^+(t) = \rho_1(\theta_{t + \frac{1}{2}})
\widetilde{\sigma}_1^- (t + \frac{1}{2}).$$ But this relation is
precisely (d1) for $\OH_1$. Therefore, the proof is concluded.
\end{proof}

\begin{remark} \label{Observacion-1-Periodicidad}
$(\textnormal{d1} \Rightarrow \textnormal{c1})$ is the only
implication that fails for $1$-periodic $\Po$-functions.
\end{remark}

\begin{theorem} \label{Teorema-Rodrigues}
Let $\rho \in \mathcal{M}(\Omega)$ for some $\Po$-invariant open
set $\Omega$ and assume that $\rho$ satisfies the functional
equation $$\textstyle \rho (\theta_t) \widetilde{\sigma}^+ (t) =
\rho(\theta_{t+1}) \widetilde{\sigma}^-(t+1).$$ Then there exists
a family of functions $\rho_j \in \mathcal{M}(\Omega)$, with $j
\ge 0$ and such that $\rho_0 = \rho$, determined by any of the
identities
\begin{eqnarray*}
\rho_{j+1}(\theta_t) & = & \textstyle \rho_j(\theta_{t -
\frac{1}{2}}) \, \widetilde{\sigma}^+ \big( t + \frac{(j-1)}{2}
\big), \\ \rho_{j+1}(\theta_t) & = & \textstyle \rho_j(\theta_{t +
\frac{1}{2}}) \, \widetilde{\sigma}^- \big( t - \frac{(j-1)}{2}
\big).
\end{eqnarray*}
Moreover, this family of functions provide the following Rodrigues
formula for the eigenfunctions $f_0, f_1, \ldots f_d$ of $\OH$
$$\Big( \prod_{j=0}^{k-1} (\mu_j - \mu_k) \Big) \rho(x) f_k (x) =
\partial_k f_k \Big( \prod_{j=0}^{k-1} \frac{\lambda^{k-j} -
\lambda^{j-k}}{\lambda - \lambda^{-1}} \Big) \Di^k \rho_k(x).$$
\end{theorem}

\begin{proof}
The recurrence relations given above for $\rho_j$ follow from
Lemmas \ref{Lema-Funciones-Auxiliares} and
\ref{Lema-Equivalencias}. Besides, by Lemma
\ref{Lema-Equivalencias}, we have $\rho_j (x) \OH_j f (x) = \Di
(\rho_{j+1} \Di f)(x)$ for any $f \in \mathcal{M}(\Omega)$. If we
combine this relation with
(\ref{Ecuacion-Autofunciones-Iteradas}), we get $$(\mu_j - \mu_k)
\rho_j (x) f_k^{(j)} (x) = \Di (\rho_{j+1} f_k^{(j+1)}) (x).$$
Then we differentiate successively this expression to obtain
$$\Big( \prod_{j=0}^{k-1} (\mu_j - \mu_k) \Big) \rho(x) f_k (x) =
\Di^k (\rho_k f_k^{(k)}) = \partial_k f_k \Big( \prod_{j=0}^{k-1}
\frac{\lambda^{k-j} - \lambda^{j-k}}{\lambda - \lambda^{-1}} \Big)
\Di^k \rho_k(x)$$ where we have applied Proposition
\ref{Proposicion-TerminoDominante}. This completes the proof.
\end{proof}

\begin{remark} \label{Observacion-Display2}
Let us recall that, assuming that $x$ is the base point of a
$\Po$-sequence with $x_{j - \frac{k}{2}}$ pairwise distinct for $j
=0,1, \ldots k$ and with the aid of Remark
\ref{Observacion-Display1}, we can display Rodrigues formula as
follows $$\Di^k \rho_k (x) = \Big( \prod_{j=0}^{k-1}
\frac{\lambda^{k-j} - \lambda^{j-k}}{\lambda - \lambda^{-1}} \Big)
\sum_{j=0}^k \frac{\rho_k(x_{j - \frac{k}{2}})}{\prod_{i \neq j}
(x_{j - \frac{k}{2}} - x_{i - \frac{k}{2}})}.$$
\end{remark}

\begin{remark}
The recurrence relations for the functions $\rho_j$ given in
Theorem \ref{Teorema-Rodrigues} have the corresponding analogs in
\cite[3.2.11 and Pg 64]{NSU}. However, the authors of \cite{NSU}
work with the functions $\rho_j^{\mbox{{\scriptsize NSU}}}$
defined by $$\rho_j^{\mbox{{\scriptsize NSU}}} (s) =
\rho_j(\theta_{s + \frac{j}{2}}).$$ Hence, these functions can not
be easily regarded as functions of the variable $\theta_s$. One of
our main contributions is that we show how to find meromorphic
solutions $\rho$ and $\rho_j$ (with the exception of the form
\textbf{C}, see Section 4) in the variable $\theta_s$. To that
aim, it is essential to regard these recurrence relations in its
most intrinsic form as in (c2) of Lemma \ref{Lema-Equivalencias}.
This result has no analog in \cite{NSU}.
\end{remark}

\begin{remark}
Another remarkable advantage of this alternative formulation is
that affine transformations become trivial. A good example to
illustrate this is Fischer's paper \cite{F}, where the author gets
a Rodrigues formula for the $q$-Racah polynomials. Fischer
combines results from \cite{AW} and \cite{NSU} with some
independent computations. However, from our point of view, these
kind of results become trivial consequences. Namely, it is very
well-known that the $q$-Racah polynomials $\mathrm{R}_k$ are
related to the Askey-Wilson polynomials $\mathrm{P}_k$ via
$$\mathrm{R}_k(\overline{x}; \alpha, \beta, \gamma, \delta) =
\mathrm{P}_k (x; a, b, c, d)$$ where $\overline{x} = 2 \sqrt{q
\gamma \delta}\, x$ and $$\begin{array}{ll} a = \sqrt{q \gamma
\delta} & \qquad b = \alpha \sqrt{q / \gamma \delta} \\ c = \beta
\sqrt{\delta q / \gamma} & \qquad d = \sqrt{\gamma q / \delta}.
\end{array}$$ Therefore, if we are given a Rodrigues formula for
the Askey-Wilson polynomials (see Section \ref{Section4}) and we
want a Rodrigues formula for the $q$-Racah polynomials, it
suffices to take $\overline{\rho}_k(\overline{x}) = \rho_k (x)$
and substitute the operator D, associated to the $\mathrm{P}_k$'s,
by the operator $\overline{\mathrm{D}}$ associated to the
$\mathrm{R}_k$'s.
\end{remark}

\begin{remark} \label{Observacion-Caso-Continuo}
The version of Theorem \ref{Teorema-Rodrigues} for the continuous
canonical form, for which every $\Po$-sequence is constant and
hence $1$-periodic, is very well-known. In this case we have
$\sigma_j(x) = \sigma(x)$ and $\tau_j(x) = \tau(x) + j
\sigma'(x)$. Moreover, the functional equations given in Theorem
\ref{Teorema-Rodrigues} must be replaced by $$\frac{d (\rho
\sigma)}{dx} (x) = \rho(x) \tau(x) \qquad \mbox{and} \qquad
\rho_j(x) = \rho(x) \sigma^j(x).$$ This also follows from Lemmas
\ref{Lema-Conmutacion} and \ref{Lema-Equivalencias}. But a direct
proof is quite simpler.
\end{remark}

\subsection{Discrete orthogonality relations}
\label{Subsection-DOR}

In this paragraph we shall assume the existence of meromorphic
functions $\rho$ and $\rho_j$ satisfying the hypothesis of Theorem
\ref{Teorema-Rodrigues} and we show how these functions can be
used to obtain discrete orthogonality relations. Given $f, g \in
\Pol[x]$, we define $h(u,v)$ to be the symmetric polynomial which
coincides with $$\frac{f(u) g(v) - f(v) g(u)}{v-u}$$ whenever $u
\neq v$. Then, we define the Wronskian $\mathrm{W}: \Pol[x] \times
\Pol[x] \rightarrow \Pol[x]$ to be the bilinear form
$$\mathrm{W}(f,g)(x) = h(A(x) + \sqrt{\delta(x)}, A(x) -
\sqrt{\delta(x)}).$$ Recall that $\mathrm{W}(f,g) \in \Pol[x]$ for
any $f, g \in \Pol[x]$ since $h(u,v)$ is symmetric in $(u,v)$.

\begin{lemma} \label{Lema-Wronski}
We have $\rho (x) (f \OH g - g \OH f)(x) = \Di (\rho_1
\mathrm{W}(f,g))(x)$ for all $x \in \Omega$.
\end{lemma}

\begin{proof}
Given a point $x \in \Omega$, we know by Lemma
\ref{Lema-Equivalencias} that $$\rho (x) (f \OH g - g \OH f)(x) =
(f \Di(\rho_1 \Di g) - g \Di(\rho_1 \Di f)) (x).$$ On the other
hand, let us consider a $\Po$-sequence with base point $x$. Then
we can assume by a continuity argument that $x_{t + \frac{1}{2}}
\neq x_{t - \frac{1}{2}}$ for $t = 0, \pm \frac{1}{2}$. In that
case, since $x = x_0$, we can rewrite the right hand side as
follows
\begin{eqnarray*}
\mathrm{RHS} & = & \frac{f(x_0)}{x_{1/2} - x_{-1/2}} \Big[
\rho_1(x_{1/2}) \frac{g(x_1) - g(x_0)}{x_1 - x_0} -
\rho_1(x_{-1/2}) \frac{g(x_0) - g(x_{-1})}{x_0 - x_{-1}} \Big]
\\ & - & \frac{g(x_0)}{x_{1/2} - x_{-1/2}} \Big[ \rho_1(x_{1/2})
\frac{f(x_1) - f(x_0)}{x_1 - x_0} - \rho_1(x_{-1/2}) \frac{f(x_0)
- f(x_{-1})}{x_0 - x_{-1}} \Big]
\end{eqnarray*}
Since the terms in $f(x_0) g(x_0)$ cancel, we easily obtain the
desired relation.
\end{proof}

\begin{theorem} \label{Teorema-Ortogonal-Discreto}
Let us consider a $\Po$-sequence $(x_j)$ and a non-negative
integer $m$. Then, if $\rho$ is regular at $x_0, x_1, \ldots, x_m$
and $\rho_1(x_{-1/2}) = \rho_1(x_{m + \frac{1}{2}}) = 0$, then the
eigenfunctions $f_0, f_1, \ldots, f_m$ satisfy the following
orthogonality relations $$\langle f_j, f_k \rangle = \sum_{i =
0}^m f_j(x_i) f_k(x_i) \rho(x_i) (x_{i + \frac{1}{2}} - x_{i -
\frac{1}{2}}) = 0 \qquad \mbox{for} \quad j \neq k.$$
\end{theorem}

\begin{proof}
Lemma \ref{Lema-Wronski} gives
\begin{eqnarray*}
\langle f, \OH g \rangle - \langle \OH f, g \rangle & = & \sum_{i
= 0}^m \Di (\rho_1 \mathrm{W}(f,g))(x_i) (x_{i + \frac{1}{2}} -
x_{i - \frac{1}{2}}) \\ & = & \sum_{i = 0}^m \rho_1 (x_{i +
\frac{1}{2}}) \mathrm{W}(f,g)(x_{i + \frac{1}{2}}) - \rho_1 (x_{i
- \frac{1}{2}}) \mathrm{W}(f,g)(x_{i - \frac{1}{2}})
\\ & = & \rho_1 (x) \mathrm{W}(f,g)(x) \Big|_{x =
x_{-\frac{1}{2}}}^{x = x_{m + \frac{1}{2}}} = 0.
\end{eqnarray*}
Taking $f = f_j$, $g = f_k$ and recalling that $\mu_j \neq \mu_k$
for $j \neq k$, the result follows.
\end{proof}

\begin{remark}
The inner product described in Theorem
\ref{Teorema-Ortogonal-Discreto} can be complex.
\end{remark}

\begin{remark} In Section \ref{Section4} we shall study how free
we are at the time of choosing the functions $\rho$ and $\rho_j$.
Using this it can be checked that, when
$$\mbox{ord}_{x_{-1/2}}(\rho_1)
> \mbox{ord}_{x_0}(\rho_1) = \mbox{ord}_{x_1}(\rho_1) = \cdots =
\mbox{ord}_{x_m}(\rho_1) < \mbox{ord}_{x_{m + 1/2}}(\rho_1),$$ we
can modify the election of $\rho$ and $\rho_1$ so that hypotheses
of Theorem \ref{Teorema-Ortogonal-Discreto} hold.
\end{remark}

\section{Explicit formulas for the canonical forms}
\label{Section4}

In this section we apply our methods to study the hypergeometric
polynomials that arise when $\Po$ takes one of the following
canonical forms: \textbf{T}, \textbf{G}, \textbf{Q}, \textbf{A} or
\textbf{C}. Although the canonical forms \textbf{O} and \textbf{E}
could also be treated with our methods, we shall not cover them
here, as we explained at the beginning of Section \ref{Section3}.
In each case we shall follow the following steps:
\begin{itemize}
\item[\bf 1.] \textbf{Parameterization} of $\sigma$ and $\tau$ in
terms of certain (Laurent) polynomial(s).
\item[\bf 2.] \textbf{Expressions} for $\alpha_2$, $\beta_1$,
$\sigma_j$, $\tau_j$ and $\mu_k$ in terms of our parameterization.
\item[\bf 3.] \textbf{Taylor formulas} for the eigenfunctions
$f_0, f_1, \ldots f_d$.
\item[\bf 4.] \textbf{Rodrigues formulas} for the eigenfunctions
$f_0, f_1, \ldots f_d$.
\item[\bf 5.] \textbf{Classical families} of hypergeometric
polynomials which belong to that form.
\end{itemize}
When dealing with Rodrigues formulas, we shall give in each case a
method to construct the function $\rho$. The iterated functions
$\rho_j$ can be easily obtained from $\rho$ by applying Theorem
\ref{Teorema-Rodrigues}. On the other hand, we shall assume again
in what follows that the eigenvalues $\mu_k$ are pairwise
distinct.

\subsection{Continuous case}

$\Po(x,y) = x^2 + y^2 - 2xy$ and $\theta_t = x_0$.

\vskip0.2cm

\noindent \textbf{1.} In this case we do not need any
parameterization.

\vskip0.2cm

\noindent \textbf{2.} It can be checked that we have $$\mu_k -
\mu_i = (i-k) (\alpha_2 (i+k) - \alpha_2 + \beta_1)$$ $$\sigma_j
(x) = \sigma (x), \quad \tau_j(x) = \tau(x) + j \sigma'(x).$$

\noindent \textbf{3.} Taking $x_0$ so that $\sigma(x_0) = 0$, we
obtain from Theorem \ref{Teorema-Expresiones-Polinomicas}
$$\partial_j f_k (x_0) = \partial_k f_k \binom{k}{j}
\prod_{i=j}^{k-1} \frac{\tau(x_0) + i \sigma'(x_0)}{\alpha_2 (k+i)
- \alpha_2 + \beta_1}.$$

\noindent \textbf{4.} As we have pointed out in Remark
\ref{Observacion-Caso-Continuo}, the Rodrigues formulas for this
case are given by the solutions of the functional equations
$$\frac{d(\rho \sigma)}{dx}(x) = \rho(x) \tau(x) \qquad \mbox{and}
\qquad \rho_j(x) = \rho(x) \sigma^j(x).$$ In general, we can not
avoid to obtain ramified solutions. Quite surprisingly, this is
the unique canonical form satisfying this property. Namely, as we
shall see below, for any other canonical form we can always find
non-ramified solutions of the corresponding functional equation
for $\rho$.

\vskip0.2cm

\noindent \textbf{5.} Taking $\sigma(x) = x^2 - 1$ (resp.
$\sigma(x) = x$ or $\sigma(x) = 1$), we obtain the Jacobi (resp.
Laguerre or Hermite) polynomials. Note that, if $\sigma(x) = 1$
for all $x \in \C$, there are no roots of $\sigma$. In other
words, the Hermite polynomials do not arise from the expressions
given above and other methods are required to obtain them.
However, the corresponding Rodrigues formula provided by Theorem
\ref{Teorema-Rodrigues} and Remark \ref{Observacion-Caso-Continuo}
also holds in this case. As we shall see in our analysis for the
\textbf{T} canonical form, the $q$-Hermite polynomials behave in
the same fashion.

\subsection{Arithmetic case}

$\Po(x,y) = (x - y + \frac{1}{2})(x - y - \frac{1}{2})$ and
$\theta_t = x_0 \pm t$.

\vskip0.2cm

\noindent \textbf{1.} Let us consider two polynomials
$\chi^{\pm}(x) = \gamma_2 x^2 + \gamma_1^{\pm} x + \gamma_0^{\pm}$
with the same coefficient for $x^2$. Clearly, there exists a
bijective correspondence between these kind of pairs and
$\Pol_2[x] \times \Pol_1[x]$ given by $$\chi^+(x) = \sigma(x) +
\frac{\tau(x)}{2} \qquad \sigma(x) = \frac{\chi^+(x) +
\chi^-(x)}{2}$$ $$\chi^-(x) = \sigma(x) - \frac{\tau(x)}{2} \qquad
\tau(x) = \chi^+(x) - \chi^-(x).$$

\noindent \textbf{2.} We have $\mu_k - \mu_i = (i-k) \big(
\gamma_2(i+k) + \gamma_1^+ - \gamma_1^- - \gamma_2 \big)$,
$\alpha_2 = \gamma_2$ and $\beta_1 = \gamma_1^+ - \gamma_1^-$.
Taking the $\Po$-function $\theta_t = t$, we easily obtain the
relations $\widetilde{\sigma}^{\pm}(t) = \chi^{\pm}(t)$. Thus, it
turns out by Lemma \ref{Lema-Funciones-Auxiliares} that
$\widetilde{\sigma}_j^{\pm}(t) = \chi^{\pm} (t \pm \frac{j}{2})$.
Therefore, $\sigma_j$ and $\tau_j$ are parameterized by
$$\textstyle \chi_j^{\pm}(x) = \chi^{\pm} (x \pm \frac{j}{2}).$$
In particular, we can write $\tau_j(x) = \chi^+(x + \frac{j}{2}) -
\chi^-(x - \frac{j}{2})$.

\vskip0.2cm

\noindent \textbf{3.} We have to consider two essentially
different kind of $\Po$-functions:
\begin{itemize}
\item[(a)] If $\theta_t = x_0 + t$, we have
$\widetilde{\sigma}^-(t) = \chi^-(x_0 + t)$. Thus, we take $x_0$
such that $\chi^-(x_0) = 0$ to satisfy the hypothesis of Theorem
\ref{Teorema-Expresiones-Polinomicas}. In that case,
$\tau_i(\theta_{i/2}) = \chi^+(x_0 + i)$ and Theorem
\ref{Teorema-Expresiones-Polinomicas} gives $$\partial_j f_k
(\theta_{j/2}) = \partial_k f_k \binom{k}{j} \prod_{i=j}^{k-1}
\frac{\chi^+(x_0 + i)}{\gamma_2(i+k) + \gamma_1^+ - \gamma_1^- -
\gamma_2}.$$
\item[(b)] If $\theta_t = x_0 - t$, we have
$\widetilde{\sigma}^-(t) = \chi^+(x_0 - t)$. Thus, if we choose
the point $x_0$ so that $\chi^+(x_0) = 0$, we obtain
$\tau_i(\theta_{i/2}) = - \chi^-(x_0 - i)$. In particular,
$$\partial_j f_k (\theta_{j/2}) = \partial_k f_k \binom{k}{j}
\prod_{i=j}^{k-1} \frac{- \chi^-(x_0 - i)}{\gamma_2(i+k) +
\gamma_1^+ - \gamma_1^- - \gamma_2}.$$
\end{itemize}
Every case with pairwise distinct eigenvalues gives rise to an
expression of this kind.

\vskip0.2cm

\noindent \textbf{4.} Taking the $\Po$-function $\theta_t = t$,
the Rodrigues formulas given by Theorem \ref{Teorema-Rodrigues}
arise from the solutions of
\begin{eqnarray*}
\rho(t) \chi^+ (t) & = & \rho(t+1) \chi^-(t+1) \\ \rho_{j+1}(t) &
= & \textstyle \rho_j(t \mp \frac{1}{2}) \chi^{\pm} (t \pm
\frac{j-1}{2}).
\end{eqnarray*}
Non-vanishing entire solutions are available whenever $\chi^{\pm}$ are
not identically zero. Moreover, those solutions can be obtained by
multiplying the suitable solutions of the following root cases
$$\hskip-7pt \begin{array}{llc|cll} \left. \begin{array}{l}
\chi^+(x) = 1 \\ \chi^-(x) = x - \xi^- \end{array} \!\! \right\} &
\displaystyle \rho(x) = \frac{1}{\Gamma(1+x-\xi^-)} & & & \left.
\begin{array}{l} \chi^+(x) = \zeta \neq 0 \\ \chi^-(x) = 1
\end{array} \!\! \right\} & \displaystyle \rho(x) = \zeta^x \\
\null \\ \left. \begin{array}{l} \chi^+(x) = x - \xi^+ \\
\chi^-(x) = 1 \end{array} \!\! \right\} & \displaystyle \rho(x) =
\frac{1}{\Gamma(1-x+\xi^+)} & & & \left. \begin{array}{l}
\chi^+(x) = 1 \\ \chi^-(x) = \zeta \neq 0 \end{array} \!\!
\right\} & \displaystyle \rho(x) = \zeta^{-x}. \end{array}$$ Other
solutions are available multiplying $\rho$ by a 1-periodic
meromorphic function.

\vskip0.2cm

\noindent \textbf{5.} Taking $\chi^{\pm}(x) = (x - \xi_0^{\pm}) (x
- \xi_1^{\pm})$, we obtain the Hahn polynomials
\begin{itemize}
\item[(a)] $\displaystyle f_k(x) = f_k(\xi_0^-) \,\ _3F_2
\left( \begin{array}{c|} -k, \ \xi_0^- - x, \ k-1 + \xi_0^- +
\xi_1^- - \xi_0^+ - \xi_1^+ \\ \xi_0^- - \xi_0^+, \ \xi_0^- -
\xi_1^+ \end{array} \ 1 \right)$
\item[(b)] $\displaystyle f_k(x) = f_k(\xi_0^+) \,\ _3F_2 \left(
\begin{array}{c|} -k, \ x - \xi_0^+, \ k-1 + \xi_0^- + \xi_1^- -
\xi_0^+ - \xi_1^+ \\ \xi_0^- - \xi_0^+, \ \xi_1^- - \xi_0^+
\end{array} \ 1 \right).$
\end{itemize}

\subsection{Quadratic case}

$\Po(x,y) = x^2 + y^2 - 2xy  - \frac{1}{2}(x + y) + \frac{1}{16}$
and $\theta_t = (t + t_0)^2$.

\vskip0.2cm

\noindent \textbf{1.} Now we consider a polynomial $\chi(t) =
\gamma_4 t^4 + \gamma_3 t^3 + \gamma_2 t^2 + \gamma_1 t +
\gamma_0$. Then, the parameterization is given by $$\chi(t) =
\sigma(t^2) + t \tau(t^2)$$ $$\sigma(t^2) = \frac{\chi(t) +
\chi(-t)}{2} \qquad \mbox{and} \qquad \tau(t^2) = \frac{\chi(t) -
\chi(-t)}{2t}.$$

\noindent \textbf{2.} We have $\alpha_2 = \gamma_4$, $\beta_1 =
\gamma_3$ and $$\mu_k - \mu_i = (i-k) \big( \gamma_4(i+k) -
\gamma_4 + \gamma_3 \big).$$ Taking $\theta_t = t^2$ and applying
Lemma \ref{Lema-Funciones-Auxiliares}, we get
$\widetilde{\sigma}_j^{\pm}(t) = \chi(\pm t + \frac{j}{2})$ for $j \ge
0$. Hence, $\sigma_j$ and $\tau_j$ are parameterized by the
polynomial $$\textstyle \chi_j(x) = \chi(x + \frac{j}{2}).$$ In
particular, we have $$\tau_j(t^2) = \frac{\chi(t + \frac{j}{2}) -
\chi(- t + \frac{j}{2})}{2t}.$$

\noindent \textbf{3.} It suffices to consider the $\Po$-functions
$\theta_t = (t + t_0)^2$. In that case, we have
$\widetilde{\sigma}^-(t) = \chi(-t - t_0)$ and so we impose
$\chi(- t_0) = 0$. This condition provides the factorization
$\chi(t) = (t + t_0) \chi_0(t)$ for some $\chi_0 \in \Pol_3[x]$.
Therefore, it turns out that $\tau_i(\theta_{i/2}) = \chi_0(t_0 +
i)$. Theorem \ref{Teorema-Expresiones-Polinomicas} gives
$$\partial_j f_k (\theta_{j/2}) = \partial_k f_k \binom{k}{j}
\prod_{i=j}^{k-1} \frac{\chi_0(t_0 + i)}{\gamma_4(i+k) - \gamma_4
+ \gamma_3}.$$ Every case with pairwise distinct eigenvalues gives
rise to an expression of this kind.

\vskip0.2cm

\noindent \textbf{4.} Taking the $\Po$-function $\theta_t = t^2$
and following Theorem \ref{Teorema-Rodrigues}, we have to solve
\begin{eqnarray*}
\rho(t^2) \chi(t) & = & \rho((t+1)^2) \chi(-t-1) \\
\rho_{j+1}(t^2) & = & \textstyle \rho_j((t \mp \frac{1}{2})^2)
\chi (\pm t + \frac{j-1}{2}).
\end{eqnarray*}
Non-vanishing entire solutions are always available. Namely, any
solution arise as a product of functions  like $$\rho(t^2) =
\frac{1}{\Gamma(1-\xi+t) \Gamma(1-\xi-t)},$$ which solves the
functional equation for $\chi(t) = t + \xi$ and $\xi \in \C$.
Moreover, other solutions appear when we multiply $\rho$ by any
even 1-periodic meromorphic function.

\vskip0.2cm

\noindent \textbf{5.} Taking $\chi(x) = \displaystyle
\prod_{\nu=0}^3 (x + \xi_{\nu})$ and $t_0 = \xi_0$, we obtain the
Wilson polynomials $$f_k(t^2) = f_k(\xi_0^2) \,\ _4F_3 \left(
\begin{array}{c|} -k, \ \xi_0 + t, \ \xi_0 - t, \ k-1 + \xi_0 +
\xi_1 + \xi_2 + \xi_3 \\ \xi_0 + \xi_1, \ \xi_0 + \xi_2, \ \xi_0 +
\xi_3 \end{array} \ 1 \right).$$

\subsection{Some remarks on a functional equation}
\label{Subsection-EF}

Finally, it remains to study the geometric and trigonometric
canonical forms. In order to give explicit Rodrigues formulas for
these particular cases, we shall need to solve functional
equations of the following type
\begin{equation} \label{Ecuacion-Funcional}
\mathrm{F}(x) = \mathrm{R}(x) \mathrm{F}(qx),
\end{equation}
where $\mathrm{R}$ stands for a non-zero rational function. An
open subset $\Omega$ of the complex plane will be called
\textbf{$q$-invariant} when $q \Omega = \Omega$. We are interested
in non-zero meromorphic solutions defined in a $q$-invariant
connected open set $\Omega$. We shall need to apply the following
basic results.

\begin{lemma} \label{Lema-Factores}
Let $\mathrm{F}_1$ and $\mathrm{F}_2$ be non-zero meromorphic
solutions of the functional equation $(\ref{Ecuacion-Funcional})$
in a $q$-invariant domain $\Omega$ with respect to the rational
functions $\mathrm{R}_1$ and $\mathrm{R}_2$ respectively. Then
$\mathrm{F}_1 \mathrm{F}_2$ and $\mathrm{F}_1 / \mathrm{F}_2$ are
non-zero meromorphic solutions of $(\ref{Ecuacion-Funcional})$ in
$\Omega$ with respect to the rational functions $\mathrm{R}_1
\mathrm{R}_2$ and $\mathrm{R}_1 / \mathrm{R}_2$ respectively.
\end{lemma}

\begin{lemma} \label{Lema-Inversion}
A function $\mathrm{F}$ solves $(\ref{Ecuacion-Funcional})$ in a
$q$-invariant domain $\Omega$ for some $q \in \C$ if and only if
$\mathrm{F}$ solves the functional equation $$\mathrm{F}(x) =
\mathrm{S}(x) \mathrm{F}(q^{-1} x) \quad \mbox{in $\Omega$ with}
\quad \mathrm{S}(x) = \frac{1}{\mathrm{R}(q^{-1}x)}.$$
\end{lemma}

\begin{lemma} \label{Lema-Unicidad}
Let us suppose that a given $q$-invariant connected open set
$\Omega$ contains the point $0$. Then, if there exist a non-zero
meromorphic solution of $(\ref{Ecuacion-Funcional})$ in $\Omega$,
it is unique up to a constant factor.
\end{lemma}

\begin{proof}
The quotient $\mathrm{F}$ of any two meromorphic solutions of the
functional equation (\ref{Ecuacion-Funcional}) satisfies
$\mathrm{F}(x) = \mathrm{F}(qx)$ in $\Omega$. Therefore, it
suffices to write this relation in terms of the Laurent series of
$\mathrm{F}$ in a neighborhood of $0$. This completes the proof.
\end{proof}

By Lemma \ref{Lema-Inversion}, the case $|q| > 1$ can be reduced
to the case $|q| < 1$. Now we provide some particular solutions
for $|q| < 1$ which will be useful in the general case. Namely,
the $q$-shifted factorials $$\begin{array}{llll} \ g_{\xi}(x) =
(\xi x;q)_{\infty} \quad & \mbox{with} \quad & \mathrm{R}(x) = 1 -
\xi x \quad & \mbox{and} \quad \Omega = \C \\ \ h_{\xi}(x) = (q /
\xi x;q)_{\infty} \quad & \mbox{with} \quad & \displaystyle
\mathrm{R}(x) = \frac{- \xi x}{1 - \xi x} \quad & \mbox{and} \quad
\Omega = \C \setminus \{0\}.
\end{array}$$ Given a non-zero rational function
$\mathrm{R}$, we can always write $\mathrm{R}$ in the form
$$\mathrm{R}(x) = \zeta x^r \prod_{k=1}^n \frac{1 - \xi_k x}{1 -
\eta_k x}$$ with $r \in \Z$, $\zeta$ a non-zero complex number and
$\xi_k, \eta_k \in \C$ for $1 \le k \le n$ . Then, using the
particular solutions $g_{\xi}$ and $h_{\xi}$, Lemma
\ref{Lema-Factores} assures that the function $$\mathrm{F} = \Big[
\frac{g_{\gamma}h_{\gamma}}{g_1h_1} \Big] \, (g_1 h_1)^r \,
\prod_{k=1}^n g_{\xi_k} / g_{\eta_k} \qquad \mbox{with} \qquad
\gamma = (-1)^r \zeta$$ solves the corresponding functional
equation. $\mathrm{F}$ is a meromorphic function in $\C \setminus
\{0\}$. In fact, if $\zeta = 1$ and $r = 0$, $\mathrm{F}$ is
meromorphic in $\C$ and we know by Lemma \ref{Lema-Unicidad} that
$\mathrm{F}$ is the unique meromorphic solution in $\C$ up to a
constant factor.

\begin{remark} \label{Observacion-No-Meromorfas}
Let $\mathrm{G}$ be meromorphic in $\C$ and biperiodic with
periods $2 \pi i$ and $\log q$. The function $\mathrm{F}$ defined
by $\mathrm{F}(\exp y) = \mathrm{G}(y)$ is meromorphic in $\C
\setminus \{0\}$ and solves the functional equation $\mathrm{F}(x)
= \mathrm{F}(qx)$. The theory of elliptic functions provides
infinitely many functions of this kind. In particular, this shows
there is no uniqueness of meromorphic solutions in $\C \setminus
\{0\}$ of the functional equation (\ref{Ecuacion-Funcional}).
\end{remark}

\begin{remark}
By the previous discussion, when $\displaystyle q^n =
\prod_{k=1}^n (\xi_k / \eta_k)$ the function $$\mathrm{E}(z) = x^n
\prod_{k=1}^n \frac{g_{\xi_k}(x) h_{\xi_k}(x)}{g_{\eta_k}(x)
h_{\eta_k}(x)} \qquad (x = e^z)$$ is elliptic with periods $2 \pi
i$ and $\log q$. Moreover, $\mathrm{E}$ has a zero at $- \log
\xi_k$ and has a pole at $- \log \eta_k$ for $1 \le k \le n$. The
rest of zeros and poles of $\mathrm{E}$ can be obtained by
periodicity. Therefore we have obtained, in our terminology, a
proof of a classical existence result. This explicit expression
allows to explore the consequences derived from the non-uniqueness
of solutions $\rho$ and $\rho_j$ of the Rodrigues formulas. We
shall give more details on this topic in the following Paragraphs
4.5 and 4.6.
\end{remark}

\begin{remark}
The functional equation (\ref{Ecuacion-Funcional}) has
non-vanishing rational solutions for $q = e^{2\pi i /n}$ whenever
the following extra hypothesis hold $$\prod_{k=0}^{n-1}
\mathrm{R}(q^k x) = 1.$$ These solutions could be used to obtain
Rodrigues formulas for such values of $q$. However, we have
decided not to include them in this paper for lack of space.
\end{remark}

\begin{remark}
Putting $f(t) = \mathrm{F}(q^t)$, (\ref{Ecuacion-Funcional}) takes
the form $f(t) = \mathrm{R}(q^t) f(t+1)$. This kind of equations
have been solved in \cite{NSU} for $|q| < 1$ in terms of the
$\Gamma_q$-function. In general, these solutions can not be
regarded as non-ramified solutions in $q^t$. This alternative
functional equation is also considered in \cite{RS}, where some
solutions $f(t) = \mathrm{F}(q^t)$ are given. However, the authors
seem not to be specially interested on solutions which are
functions of the variable $q^t$ and they do not observe the role
of elliptic functions here. For more on the connection between
$q$-calculus and elliptic function theory, we refer the reader to
the references \cite{AAR} and \cite{BrI}.
\end{remark}

\subsection{Geometric case}

$\Po(x,y) = x^2 + y^2 - 2axy$ and $\theta_t = q^{\pm t} x_0$.

\vskip0.2cm

\noindent \textbf{1.} Let us consider two polynomials $\chi^{\pm}(x) =
\gamma_2^{\pm} x^2 + \gamma_1^{\pm} x + \gamma_0$ with the same value at
$x=0$. Then, the parameterization is given by $$\chi^+(x) =
\sigma(x) + \frac{\lambda - \lambda^{-1}}{2} \, x \, \tau(x)
\qquad \sigma(x) = \frac{\chi^+(x) + \chi^-(x)}{2}$$ $$\chi^-(x) =
\sigma(x) - \frac{\lambda - \lambda^{-1}}{2} \, x \, \tau(x)
\qquad \tau(x) = \frac{\chi^+(x) - \chi^-(x)}{(\lambda -
\lambda^{-1}) x}.$$

\noindent \textbf{2.} We have $\alpha_2 = (\gamma_2^+ +
\gamma_2^-)/2$, $\beta_1 = (\gamma_2^+ - \gamma_2^-)/(\lambda -
\lambda^{-1})$ and
$$\mu_k - \mu_i = \frac{\lambda^{i-k} - \lambda^{k-i}}{(\lambda -
\lambda^{-1})^2} (\gamma_2^+ \lambda^{i+k-1} - \gamma_2^-
\lambda^{1-i-k}).$$ Taking $\theta_t = q^t$, it turns out that
$\widetilde{\sigma}^{\pm}(t) = \chi^{\pm}(q^t)$ and
$\widetilde{\sigma}_j^{\pm}(t) = \chi^{\pm}(\lambda^{\pm j} q^t)$.
In other words, the polynomials $\sigma_j$ and $\tau_j$ are
parameterized by $$\chi_j^{\pm}(x) = \chi^{\pm} (\lambda^{\pm j}
x).$$ In particular, we obtain $$\tau_j(x) =
\frac{\chi^+(\lambda^j x) - \chi^-(\lambda^{-j} x)}{(\lambda -
\lambda^{-1}) x}.$$

\noindent \textbf{3.} We have to consider two essentially
different kind of $\Po$-functions:
\begin{itemize}
\item[(a)] If $\theta_t = q^t x_0$, we have
$\widetilde{\sigma}^-(t) = \chi^-(q^t x_0)$. Thus, we require
$\chi^-(x_0) = 0$ to satisfy the hypothesis of Theorem
\ref{Teorema-Expresiones-Polinomicas}. Under this assumption, we
have $\tau_i(\theta_{i/2}) = \chi^+(q^i x_0) / (\lambda -
\lambda^{-1}) \lambda^i x_0$ and Theorem
\ref{Teorema-Expresiones-Polinomicas} gives $$\partial_j f_k
(\theta_{j/2}) = \partial_k f_k \Big[ \! \begin{array}{c} k \\ j
\end{array} \! \Big]_q q^{- 2 [\binom{k}{2} - \binom{j}{2}]}
\prod_{i=j}^{k-1} \frac{\chi^+(q^i x_0) / x_0}{\gamma_2^+ -
\gamma_2^- q^{1-i-k}}.$$
\item[(b)] If $\theta_t = q^{-t} x_0$, we have
$\widetilde{\sigma}^-(t) = \chi^+(q^{- t} x_0)$. Thus, we take
$x_0$ such that $\chi^+(x_0) = 0$. This gives
$\tau_i(\theta_{i/2}) = - \lambda^i \chi^-(q^{-i} x_0) / (\lambda
- \lambda^{-1}) x_0$. In particular, $$\partial_j f_k
(\theta_{j/2}) = \partial_k f_k \Big[ \! \begin{array}{c} k \\ j
\end{array} \! \Big]_q q^{- [\binom{k}{2} - \binom{j}{2}]}
\prod_{i=j}^{k-1} \frac{- \chi^-(q^{-i} x_0) / x_0}{\gamma_2^+ -
\gamma_2^- q^{1-i-k}}.$$
\end{itemize}
Every case with pairwise distinct eigenvalues gives rise to an
expression of this kind. When $x_0 = 0$, the terms
$\chi^{\pm}(q^{\pm i}x_0) / x_0$ must be replaced by the obvious
limits.

\vskip0.2cm

\noindent \textbf{4.} Taking the $\Po$-function $\theta_t = q^t$
and following Theorem \ref{Teorema-Rodrigues}, we have to obtain
non-zero meromorphic solutions of
\begin{eqnarray*}
\rho(x) \chi^+(x) & = & \rho(qx) \chi^-(qx) \\ \rho_{j+1}(x) & = &
\rho_j(\lambda^{\mp 1} x) \chi^{\pm} (\lambda^{\pm (j-1)} x).
\end{eqnarray*}
in a $\Po$-invariant domain $\Omega$. Non-vanishing meromorphic
solutions in $\C \setminus \{0\}$ are available whenever $|q| < 1$
and $\chi^{\pm}$ are not identically zero. It suffices to multiply the
solutions of the following root cases $$\begin{array}{ll} \left.
\begin{array}{l} \chi^+(x) = 1 - \xi^+ x \\ \chi^-(x) = 1 - \xi^- x
\end{array} \!\! \right\} & \displaystyle \rho(x) = \frac{(\xi^- q
x; q)_{\infty}}{(\xi^+ x; q)_{\infty}}, \\ \null \\ \left.
\begin{array}{l} \chi^+(x) = x^{m_+} \zeta^+ \\ \chi^-(x) =
x^{m_-} \zeta^-
\end{array} \right\} & \displaystyle \rho(x) = \frac{g_{\gamma}
h_{\gamma}}{g_1h_1} (g_1h_1)^m,\end{array}$$ where $m = m_- - m_+$
and $\gamma = (-1)^m \zeta^- q^{m_-} / \zeta^+$. Besides, Remark
\ref{Observacion-No-Meromorfas} provides infinitely many
meromorphic solutions in $\C \setminus \{0\}$. For particular
cases, choosing the simplest solution might not be completely
trivial.

\vskip0.2cm

\noindent \textbf{5.} Taking $\chi^{\pm}(x) = (1 - \xi_0^{\pm} x)
(1 - \xi_1^{\pm} x)$ with $\xi_0^{\pm} \neq 0$ and $\xi_1^{\pm}
\neq 0$, we obtain the $q$-Hahn polynomials
\begin{itemize}
\item[(a)] $\displaystyle f_k(x) = f_k(1/\xi_0^-) \,\ _3\phi_2
\left( \begin{array}{c|} q^{-k}, \ 1 / \xi_0^- x, \ q^{k-1}
\xi_0^+ \xi_1^+ / \xi_0^- \xi_1^- \\ \xi_0^+ / \xi_0^-, \ \xi_1^+
/ \xi_0^- \end{array} \ q, \ q \xi_1^- x \right)$.
\item[(b)] $\displaystyle f_k(x) = f_k(1/\xi_0^+) \,\ _3\phi_2
\left( \begin{array}{c|} q^{-k}, \ \xi_0^+ x, \ q^{k-1} \xi_0^+
\xi_1^+ / \xi_0^- \xi_1^- \\ \xi_0^+ / \xi_0^-, \ \xi_0^+ /
\xi_1^- \end{array} \ q, \ q \right).$
\end{itemize}

\subsection{Trigonometric case}

$\Po(x,y) = x^2 + y^2 - 2axy  + a^2 - 1$ and $\theta_t =
\displaystyle \frac{q^t u_0 + q^{-t} u_0^{-1}}{2}$.

\vskip0.2cm

\noindent \textbf{1.} In this case, the Laurent polynomial
$\chi(u) = \gamma_{-2} u^{-2} + \gamma_{-1} u^{-1} + \gamma_0 +
\gamma_1 u + \gamma_2 u^2$ provides the following parameterization
$$\chi(u) = \sigma \Big( \frac{u + u^{-1}}{2} \Big) +
\frac{(\lambda - \lambda^{-1}) (u - u^{-1})}{4} \tau \Big( \frac{u
+ u^{-1}}{2} \Big)$$ $$\sigma \Big(\frac{u + u^{-1}}{2} \Big) =
\frac{\chi(u) + \chi(u^{-1})}{2} \quad \mbox{and} \quad \tau \Big(
\frac{u + u^{-1}}{2} \Big) = 2 \frac{\chi(u) -
\chi(u^{-1})}{(\lambda - \lambda^{-1})(u - u^{-1})}.$$

\noindent \textbf{2.} We have $\alpha_2 = 2 (\gamma_2 +
\gamma_{-2})$, $\beta_1 = 4 (\gamma_2 - \gamma_{-2}) / (\lambda -
\lambda^{-1})$ and $$\mu_k - \mu_i = 4 \frac{\lambda^{i-k} -
\lambda^{k-i}}{(\lambda - \lambda^{-1})^2} (\gamma_2
\lambda^{i+k-1} - \gamma_{-2} \lambda^{1-i-k}).$$ Taking $\theta_t
= (q^t + q^{-t}) / 2$, we obtain $\widetilde{\sigma}^{\pm}(t) =
\chi(q^{\pm t})$ and $\widetilde{\sigma}_j^{\pm}(t) =
\chi(\lambda^j q^{\pm t})$. That is, the polynomials $\sigma_j$
and $\tau_j$ are parameterized by $$\chi_j(u) = \chi(\lambda^j
u).$$ In particular, we have $$\tau_j \Big( \frac{u + u^{-1}}{2}
\Big) = 2 \frac{\chi(\lambda^j u) - \chi(\lambda^j
u^{-1})}{(\lambda - \lambda^{-1})(u - u^{-1})}.$$

\noindent \textbf{3.} It suffices to consider the $\Po$-functions
$\theta_t = (q^t u_0 + q^{-t} u_0^{-1}) / 2$. In that case, we
have $\widetilde{\sigma}^-(t) = \chi(q^{-t} u_0^{-1})$. Thus, we
take $u_0$ so that $\chi(u_0^{-1}) = 0$. This condition provides
the factorization $\chi(u) = (u - u_0^{-1}) u^{-2} \chi_0(u)$ for
some $\chi_0 \in \Pol_3[x]$. Therefore, it turns out that
$$\tau_i(\theta_{i/2}) = \frac{2 \chi_0(q^i u_0)}{\lambda^{3i}
u_0^2 (\lambda - \lambda^{-1})}.$$ Theorem
\ref{Teorema-Expresiones-Polinomicas} gives $$\partial_j f_k
(\theta_{j/2}) = \partial_k f_k \Big[ \!
\begin{array}{c} k \\ j \end{array} \! \Big]_q q^{-3 [\binom{k}{2} -
\binom{j}{2}]} \prod_{i=j}^{k-1} \frac{\chi_0(q^i u_0) /
u_0^2}{2(\gamma_2 - \gamma_{-2} q^{1-i-k})}.$$ Let us notice that,
for $\chi(u) = \gamma_{\pm 2} u^{\pm 2}$, our methods provide a
Rodrigues formula. However there is not any $u_0$ such that
$\chi(u_0^{-1}) = 0$. Hence we can not apply Theorem
\ref{Teorema-Expresiones-Polinomicas}. These cases correspond to
the $q$-Hermite polynomials.

\vskip0.2cm

\noindent \textbf{4.} Now we apply the ideas of Paragraph
\ref{Subsection-EF} to discuss some Rodrigues formulas for the
trigonometric case when $|q| < 1$. We begin by recalling that the
relations $$\Omega = \Big\{ \frac{u + u^{-1}}{2}: \ u \in \Omega'
\Big\} \qquad \mbox{and} \qquad \Omega' = \Big\{ u: \ \frac{u +
u^{-1}}{2} \in \Omega \Big\}$$ provide a bijective correspondence
between the class of $\Po$-invariant open sets $\Omega$ of the
complex plane and the class of $\lambda$-invariant open subsets
$\Omega'$ of $\C \setminus \{0\}$ satisfying that $u \in \Omega'$
if and only if $u^{-1} \in \Omega'$. Moreover, given $\Omega$ and
$\Omega'$ as above, the relation $$\mathrm{F}(u) = \mathrm{G}
\Big( \frac{u + u^{-1}}{2} \Big)$$ establishes a bijection between
meromorphic functions $\mathrm{G}$ in $\Omega$ and meromorphic
functions $\mathrm{F}$ in $\Omega'$ satisfying $\mathrm{F}(u) =
\mathrm{F}(u^{-1})$. Consequently, we shall work with functions
$\rho_j$ of the form $$\rho_j \Big( \frac{u + u^{-1}}{2} \Big) =
h_j(u)$$ where $h_j$ is meromorphic in $\Omega'$ and $h_j(u) =
h_j(u^{-1})$ for all $u \in \Omega'$. Taking the $\Po$-function
$\theta_t = (q^t + q^{-t}) / 2$ with $u = q^t$, we have to solve
\begin{eqnarray*}
h(u) \chi(u) & = & h(qu) \chi(q^{-1}u^{-1}) \\ h_{j+1}(u) & = &
h_j(\lambda^{\mp 1} u) \chi(\lambda^{j-1} u^{\pm 1}).
\end{eqnarray*}
We can always write the Laurent polynomial $\chi$ in the form
$$\chi(u) = \zeta u^{-s} \prod_{\nu=0}^{r-1} (1 - \xi_{\nu}u)$$
where $r-s = \max \{k: \gamma_k \neq 0\}$ and $\zeta =
\gamma_{k_0}$ such that $k_0 = \min \{k: \gamma_k \neq 0\}$. Using
the arguments of Paragraph \ref{Subsection-EF} and assuming $|q| <
1$, we observe that the functions $$h_j^{\xi}(u) = \frac{1}{(\xi
\lambda^j u ;q)_{\infty} (\xi \lambda^j / u ;q)_{\infty}}$$ are
meromorphic solutions in $\C \setminus \{0\}$ for $\chi(u) = 1 -
\xi u$ and invariant under the change $u \mapsto u^{-1}$. On the
other hand, when $\chi(u) = u^{-1}$ the system above reduces to
\begin{eqnarray*}
h(u) & = & q u^2 h(qu) \\ h_{j+1}(u) & = & \lambda^{1-j} u^{\pm 1}
h_j(\lambda^{\pm 1} u).
\end{eqnarray*}
Let $h^{\pm}(u) = (\mp \lambda^{1/2} u ;\lambda)_{\infty} (\mp
\lambda^{1/2} / u ;\lambda)_{\infty}$. We have $h^{\pm}(u) = \pm
\lambda^{1/2} u h^{\pm}(\lambda u)$. In particular, the functions
$h^+$ and $h^-$ solve the functional equation $h(u) = q u^2
h(qu)$. Moreover, if we define recursively the constants $k_j^{\pm}$
by $$k_{j+1}^{\pm} = \pm \lambda^{\frac{1}{2} - j} k_j^{\pm} \qquad
\mbox{with} \qquad k_0^{\pm} = 1,$$ it can be easily checked that the
functions $h_j^{\pm} = k_j^{\pm} h^{\pm}$ are meromorphic solutions in $\C
\setminus \{0\}$ for $\chi(u) = u^{-1}$ and invariant under the
change $u \mapsto u^{-1}$. In summary, by Lemma
\ref{Lema-Factores} we have meromorphic solutions in $\C \setminus
\{0\}$ of the system above for any Laurent polynomial $\chi$
$$h_j(u) = (k_j^{\pm} h^{\pm})^s (u) \prod_{\nu=0}^{r-1}
h_j^{\xi_{\nu}}(u) = (\pm 1)^{js} \lambda^{s [\frac{j}{2} -
\binom{j}{2}]} (h^{\pm})^s (u) \prod_{\nu=0}^{r-1}
h_j^{\xi_{\nu}}(u).$$ Therefore, the functions $\rho_j$ obtained
by this process are meromorphic in the complex plane. For
instance, if we have $s = 2$ and $r = 4$, we can choose the
solution $h_j^+ h_j^-$ for $\chi(u) = u^{-2}$ and so we get
\begin{equation} \label{Ecuacion-Rodriges-Askey-Wilson}
\rho_j \Big( \frac{u + u^{-1}}{2} \Big) = (-\lambda)^j q^{-
\binom{j}{2}} \frac{(\lambda u^2; q)_{\infty} (\lambda / u^2;
q)_{\infty}}{\prod_{\nu=0}^3 (\xi_{\nu} \lambda^j u ;q)_{\infty}
(\xi_{\nu} \lambda^j / u ;q)_{\infty}}.
\end{equation}
As in Remark \ref{Observacion-No-Meromorfas}, we can find other
meromorphic solutions in $\C \setminus \{0\}$ of the given
functional equation. Namely, it suffices to multiply $h$ by a
function $f$ satisfying $f(\exp x) = g(x)$ with $g$ even, $2 \pi
i$-periodic, $\log q$-periodic and meromorphic in $\C$.

\begin{remark}
The solutions we use for the functional equation $h(u) = q u^2
h(qu)$ do not follow the general solution provided in Paragraph
\ref{Subsection-EF}. The main motivation for our election has been
the simplicity of the resulting functions $\rho_k$.
\end{remark}

\noindent \textbf{5.} If $\gamma_{-2}$ and $\gamma_2$ are
non-zero, we obtain the Askey-Wilson polynomials. Namely,
normalizing if necessary, we can write $$\chi(u) = \displaystyle
u^{-2} \prod_{\nu=0}^3 (1 - \xi_{\nu} u)$$ in such a way that
$\gamma_{-2} = 1$ and $\gamma_2 = \xi_0 \xi_1 \xi_2 \xi_3$.
Therefore we get $$f_k \Big( \frac{u + u^{-1}}{2} \Big) = f_k
\Big( \frac{\xi_0 + \xi_0^{-1}}{2} \Big)\,\ _4\phi_3 \left(
\begin{array}{c|} q^{-k}, \ \xi_0 u, \ \xi_0 / u, \ \xi_0 \xi_1
\xi_2 \xi_3 q^{k-1} \\ \xi_0 \xi_1, \ \xi_0 \xi_2, \ \xi_0 \xi_3
\end{array} \ q, \ q \right).$$

\section{Concluding remarks}
\label{Section5}

\noindent {\bf 1.} In the canonical forms analyzed in Section
\ref{Section4}, the eigenfunctions are families of polynomials
depending of five parameters which correspond to the coefficients
of $\sigma$ and $\tau$. The number of free parameters can be
reduced by considering two proportional hypergeometric operators
to be equivalent. Moreover, recalling how $\OH$ is transformed
under the action of $\mbox{Aff}(\C)$ $$g \ \OH^{\Po} g^{-1} =
\zeta^2 \sigma^g (\Di^{g \cdot \Po})^2 + \zeta \tau^g \Su^{g \cdot
\Po} \Di^{g \cdot \Po},$$ we can use the isotropy subgroups of the
$\mathcal{A}_1$-orbits provided in Table I to identify different
hypergeometric operators by suitable affine transformations. In
particular, we can reduce the number of essential parameters in
each canonical form. Namely, we have two free parameters in the
continuous canonical form, three parameters in the arithmetic and
geometric cases and four parameters in the quadratic and
trigonometric canonical forms.

\vskip 5pt

\noindent {\bf 2.} The Rodrigues formula given in \cite{AW} for
the Askey-Wilson polynomials arise from
(\ref{Ecuacion-Rodriges-Askey-Wilson}) by taking $$\hskip20pt
\rho_j^{\mbox{{\scriptsize AW}}}(u) = k(u) \, \rho_j \Big( \frac{u
+ u^{-1}}{2} \Big).$$ Here, the function $k$ is defined as follows
$$k(u) = \frac{(q u^2;q)_{\infty} (q/u^2;q)_{\infty} (u -
u^{-1})}{(\lambda u^2;q)_{\infty} (\lambda/u^2;q)_{\infty}}.$$ The
function $k$ is meromorphic in $\C \setminus \{0\}$ and satisfies
the functional equation $k(qu) = k(u)$. However, we have $k(u) = -
k(u^{-1})$. Therefore, the functions used in \cite{AW} are
ramified in the variable $$x = \frac{1}{2} (u + u^{-1}).$$ The
function $k$ can be easily regarded in terms of Jacobi's $sn$
elliptic function.

\vskip 5pt

\noindent {\bf 3.} Let us consider a function $\varepsilon$
meromorphic in $\C \setminus \{0\}$ and such that
\begin{itemize}
\item[$\circ$] The function $\varepsilon$ satisfies
$\varepsilon(u) = \varepsilon(qu) = - \varepsilon(1/u)$.
\item[$\circ$] The function $\varepsilon(\lambda u) \rho_1(\frac12
(u + u^{-1}))$ is regular for $\lambda \le |u| \le \lambda^{-1}$.
\end{itemize}
Then it can be proved that, assuming that $q, \xi_{\nu} \in (0,1)$
and adapting the arguments of \cite[3.10.5.2]{NSU}, the
hypergeometric operator $\OH$ is symmetric with respect to the
product $$\langle f,g \rangle = \int_{-1}^1 f(x) g(x) \rho(x)
\widetilde{\varepsilon}(x) dx$$ where
$\widetilde{\varepsilon}(\cos (s)) = \varepsilon(e^{is})$ for $0
\le s \le \pi$. This gives the orthogonality relations for the
Askey-Wilson polynomials $f_k$. The function $k$ defined on the
previous point satisfies the required conditions on $\varepsilon$.
In fact, by the Riemann-Roch theorem, it can be checked the it is
unique up to a constant factor. Hence, this allows to recover the
well-known orthogonality relations given in \cite{AW}.

\section*{Acknowledgments}

The authors wish to thank Gabino Gonz\'{a}lez for a useful
conversation on the Riemann-Roch theorem. We also thank the
referee for some interesting comments and for having brought some
references to our attention. After this paper was accepted for
publication, Paul Terwilliger communicated to us the close
relation existing between this paper and the notion of a Leonard
pair. We refer the reader to Terwilliger's paper \cite{T} for more
on this topic. This research was supported in part by MCYT Spain
via the Project BFM 2001/0189

\enlargethispage{10pt}

\bibliographystyle{amsplain}

\end{document}